\documentclass[reqno]{amsart}

\def\draftdate{June 12, 2006}

\usepackage{amssymb}

\usepackage{xy}
\xyoption{matrix}
\xyoption{arrow}
\xyoption{curve}
\xyoption{cmtip}
\SelectTips{cm}{}
\SilentMatrices
\usepackage{ifthen}
\usepackage{xspace}

\DeclareMathAlphabet{\scr}{U}{rsfs}{m}{n}

\newcommand{\lrarrow}[1][\;]{\xrightarrow{#1}}
\newcommand{\Bmap}[4][{\lrarrow}]{\ensuremath{{#2}\colon{#3}{#1}{#4}}}
\newcommand{\map}[3]{\Bmap{#1}{#2}{#3}}
\newcommand{\vmap}[3][{\;}]{\ensuremath{{#2}{\lrarrow[{#1}]}{#3}}}
\newcommand{\backisoabove}[2][{}]{\ar@<-.5ex>@{.>}[#1]\ar@<.5ex>@{<-}[#1]^{#2}}
\newcommand{\backisobelow}[2][{}]{\ar@<.5ex>@{.>}[#1]\ar@<-.5ex>@{<-}[#1]_{#2}}

\newtheorem{theorem}{Theorem}[section]
\newtheorem{lemma}[theorem]{Lemma}
\newtheorem{corollary}[theorem]{Corollary}
\newtheorem{proposition}[theorem]{Proposition}

\theoremstyle{definition}
\newtheorem{definition}[theorem]{Definition}
\newcommand{\defterm}[1]{\emph{#1}}

\numberwithin{equation}{section}
\makeatletter
\let\c@equation\c@theorem
\makeatother

\makeatletter
\newcounter{thmpart}

\newenvironment{thmlist}[1][0]%
{\begin{list}%
{\ifthenelse{\boolean{@newlist} \and \(\not {\equal{#1}{0}}\)}%
{\hspace*{-2em} \normalfont \the\thm@headfont \upshape (\alph{thmpart})}%
{\normalfont \the\thm@headfont \upshape (\alph{thmpart})}%
}%
{\usecounter{thmpart}%
\setlength{\labelsep}{0.5em}\setlength{\leftmargin}{0pt}%
\setlength{\labelwidth}{0pt}\setlength{\itemindent}{1.8em}%
\setlength{\listparindent}{1em}}}%
{\end{list}}

\newcounter{hyppart}

\newenvironment{hyplist}[1][0]%
{\begin{list}%
{\ifthenelse{\boolean{@newlist} \and \(\not {\equal{#1}{0}}\)}%
{\hspace*{-2em} \normalfont \the\thm@headfont \upshape (\roman{hyppart})}%
{\normalfont \the\thm@headfont \upshape (\roman{hyppart})}%
}%
{\usecounter{hyppart}%
\setlength{\labelsep}{0.5em}\setlength{\leftmargin}{0pt}%
\setlength{\labelwidth}{0pt}\setlength{\itemindent}{1.8em}%
\setlength{\listparindent}{1em}}}%
{\end{list}}

\def\llabel#1{\@bsphack
  \protected@write\@auxout{}%
         {\string\newlabel{#1}{{\alph{\@listctr}}{\thepage}}}%
  \@esphack}
\makeatother

\newcommand{\PPa}{\textup{\textbf{(PP)}}\xspace}
\newcommand{\tPPa}{\textup{\textbf{(PTP)}}\xspace}
\newcommand{\SMa}{\textup{\textbf{(SM7)}}\xspace}
\newcommand{\cSMa}{\textup{\textbf{(Cot)}}\xspace}

\newcommand{\ENa}{\textup{\textbf{(Enr)}}\xspace}
\newcommand{\UF}{\textup{\textbf{(Unit)}}\xspace}
\newcommand{\UC}{$\textup{\textbf{(Unit}}'\textup{\textbf{)}}$\xspace}

\newcommand{\ENCa}{\textup{\textbf{(Enr)}}\xspace}
\newcommand{\UE}{\textup{\textbf{(HoUnit)}}\xspace}
\newcommand{\cUF}{\textup{\textbf{(CUnit)}}\xspace}
\newcommand{\tUC}{$\textup{\textbf{(TUnit)}}$\xspace}

\newcommand{\bindmp}{\mathbin{\widehat \monprod}}

\newcommand{\tL}{\mathbf{L}}
\newcommand{\tR}{\mathbf{R}}

\newcommand{\iso}{\cong}     
\newcommand{\overto}[1]{\xrightarrow{#1}}
\newcommand{\op}{\mathrm{op}}
\newcommand{\Ho}{\mathrm{Ho}}
\newcommand{\tand}{\qquad\text{and}\qquad}

\def\quickop#1{\expandafter\DeclareMathOperator\csname #1\endcsname{#1}}
\quickop{Ext} \quickop{Tor}\quickop{id}\quickop{Id}

\newcommand{\Moncat}{\ensuremath{\scr{M}}\xspace}
\newcommand{\Moncatb}{\ensuremath{\scr{N}}\xspace}
\newcommand{\Leftmod}[1][{\monoid}]{\ensuremath{{{}_{#1}\Moncat}}\xspace}
\newcommand{\Rightmod}[1][{\monoid}]{{\ensuremath{\Moncat\!{}_{#1}}}\xspace}
\newcommand{\Bimod}[2]{\ensuremath{{{}_{#1}\Moncat\!{}_{#2}}}\xspace}
\newcommand{\EC}{\ensuremath{\scr{C}}\xspace}
\newcommand{\ED}{\ensuremath{\scr{D}}\xspace}
\newcommand{\EE}{\ensuremath{\scr{E}}\xspace}
\newcommand{\EH}{\scr{H}}

\newcommand{\monprod}{\ensuremath{\wedge}}
\newcommand{\unit}{\ensuremath{I}\xspace}
\newcommand{\monoid}{\ensuremath{A}\xspace}
\newcommand{\monoidone}{\ensuremath{A}\xspace}
\newcommand{\monoidtwo}{\ensuremath{B}\xspace}
\newcommand{\monoidthree}{\ensuremath{C}\xspace}
\newcommand{\monoidfour}{\ensuremath{D}\xspace}

\newcommand{\cofunit}{I_{c}}
\newcommand{\cuf}{\tilde\ell_{c}}
\newcommand{\cum}{\ell_{c}}

\DeclareMathSymbol{\boxbin}{\mathbin}{AMSa}{"03}
\let\boxmonprod\boxbin
\newcommand{\dermonprod}[2]{{#1}\mathbin{\widehat \monprod}{#2}}
\newcommand{\emptydermonprod}{\ensuremath{\mathop{\widehat \monprod}}\xspace}

\newcommand{\func}[2]{[{#1},{#2}]}
\newcommand{\emptyfunc}{\ensuremath{\func{-}{-}}\xspace}
\newcommand{\boxfunc}[2]{{}_{\boxbin}\func{#1}{#2}}
\newcommand{\derfunc}[2]{[\![{#1},{#2}]\!]}
\newcommand{\emptyderfunc}{\ensuremath{\derfunc{-}{-}}\xspace}

\newcommand{\balprod}[1][{\monoid}]{\wedge_{#1}}

\newcommand{\WTor}{\Tor}
\newcommand{\WowTor}[5]{\WTor_{#2}({}_{#1}{#4},{#5}{}_{#3})}

\newcommand{\lmfunc}[3][{\monoid}]{{{}^{#1}[{#2},{#3}]}}
\newcommand{\emptylmfunc}[1][{\monoid}]{\ensuremath{\lmfunc[#1]{-}{-}}\xspace}
\newcommand{\rmfunc}[3][{\monoid}]{{[{#2},{#3}]{}^{#1}}}
\newcommand{\bmfunc}[4]{{{}^{#1}[{#3},{#4}]{}^{#2}}}
\newcommand{\derlmfunc}[3][{\monoid}]{{}^{#1}[\![{#2},{#3}]\!]}

\newcommand{\derbmfunc}[4]{{}^{#1}[\![{#3},{#4}]\!]{}^{#2}}
\newcommand{\emptyderlmfunc}[1][{\monoid}]{\ensuremath{\derlmfunc[#1]{-}{-}}\xspace}
\newcommand{\WExt}{\Ext}
\newcommand{\WowExt}[5]{\WExt_{#1}({#4}{}_{#2},{#5}{}_{#3})}
\newcommand{\WowExtP}[5]{\WExt_{{#3}^{\op}}({}_{#2}{#4},{}_{#1}{#5})}
\newcommand{\WowExtp}[5]{\WExt_{{#3}}({}_{#2}{#4},{}_{#1}{#5})}
\newcommand{\WowExtR}[5]{\WExt_{#1}({}_{#2}{#4},{#5}{}_{#3})}

\newcommand{\derfree}{\mathbb{F}}
\newcommand{\dercofree}{\mathbb{F}^{\sharp}}

\newcommand{\ob}[1]{\mathrm{Ob}(#1)}
\newcommand{\objectone}{C}
\newcommand{\objecttwo}{D}
\newcommand{\objectthree}{E}

\newcommand{\efunc}[3][{\EC}]{{}^{#1}[{#2},{#3}]}
\newcommand{\emptyefunc}[1][{\EC}]{\ensuremath{\efunc[{#1}]{-}{-}}\xspace}
\newcommand{\derefunc}[3][{\EC}]{{}^{#1}[\![{#2},{#3}]\!]}
\newcommand{\emptyderefunc}[1][{\EC}]{\ensuremath{\derefunc[#1]{-}{-}}\xspace}

\newcommand{\tensor}[3][{\EC}]{{#2}\otimes {#3}}
\newcommand{\emptytensor}[1][{\EC}]{\otimes}
\newcommand{\dertensor}[3][{\EC}]{{#2}\mathbin{\widehat \otimes}{#3}}
\newcommand{\emptydertensor}[1][{\EC}]{{\widehat \otimes}}
\newcommand{\cotensor}[3][{\EC}]{[{#2},{#3}]}
\newcommand{\emptycotensor}[1][{\EC}]{[{-},{-}]}
\newcommand{\dercotensor}[3][{\EC}]{[\![{#2},{#3}]\!]}
\newcommand{\emptydercotensor}[1][\EC]{\dercotensor[#1]{-}{-}}

\newcommand{\ten}[3][{\EC}]{\ensuremath{\tensor[{#1}]{#3}{#2}}\xspace}

\newcommand{\coten}[3][{\EC}]{\ensuremath{\cotensor[{#1}]{#2}{#3}}\xspace}

\newcommand{\boxten}[2]{\ensuremath{{#2}\boxtimes{#1}}\xspace}
\newcommand{\boxcot}[2]{\ensuremath{\boxfunc{#1}{#2}}\xspace}
\newcommand{\boxefunc}[3][{\EC}]{\ensuremath{{}^{#1}_{\boxbin}\func{#2}{#3}}\xspace}

\begin{document}
\title{Modules in Monoidal Model Categories}

\author{L. Gaunce Lewis, Jr.}
\address{Department of Mathematics, Syracuse University, 
Syracuse, NY \ 13244-1150}
\email{lglewis@syr.edu}

\author{Michael A. Mandell}
\address{Department of Mathematics, Indiana University,
Bloomington, IN \ 47405}
\email{mmandell@indiana.edu}
\thanks{The second author was supported in part by NSF grant DMS-0504069}

\date{\draftdate}

\begin{abstract}
This paper studies the existence of and compatibility between derived
change of ring, balanced product, and function module derived functors
on module categories in monoidal model categories.   
\end{abstract}

\maketitle


\section{Introduction}\label{IntSec}

A ``monoidal model category'' \Moncat is a category having a symmetric
monoidal closed structure and a closed model category structure which
satisfy certain compatibility conditions (reviewed in
Section~\ref{MndlModCatSec} below).  These conditions ensure that the
homotopy category $\Ho\Moncat$ inherits a symmetric monoidal structure
and that the localization functor is (lax) symmetric monoidal.
Schwede and Shipley began the study of monoidal model categories in
\cite{ssmonoidal}.  There, they provide good criteria for the
categories of modules and algebras over a monoid $\monoid$ in
$\Moncat$ to inherit closed model structures from $\Moncat$.  The
purpose of this paper is to study the existence and behavior of the
derived functors of certain commonly used functors relating various
categories of modules over a monoid in a monoidal model category.
These functors are all variants of the ``function object'' and
``balanced product'' constructions.

Let \monprod\ denote the symmetric monoidal product, \unit\ denote the
unit, and \func{-}{-} denote the internal function object for \Moncat.
A left $\monoid$-module is an object $M$ of $\Moncat$ with an
associative and unital left action map $\monoid\monprod M\to M$.  For
left $\monoid$-modules $L$ and $M$, the left $\monoid$-module function
object $\lmfunc{L}{M}$ is the equalizer in \Moncat
\[
\xymatrix{ \strut \lmfunc{L}{M} \ar[r] &\func{L}{M} \ar@<-.5ex>[r]
\ar@<.5ex>[r] &\func{\monoid \monprod L}{M}, }
\]
where one of the righthand arrows is induced by the \monoid-action
$\monoid \monprod L\to L$ and the other the composite of the map
$\func{L}{M}\to \func{\monoid \monprod L}{\monoid \monprod M}$ and the
$\monoid$-action $\monoid \monprod M\to M$.  Similarly, for a right
\monoid-module $M$ and a left $\monoid$-module $N$, the balanced
product $M\balprod N$ is the coequalizer in \Moncat
\[
\xymatrix{ M \monprod \monoid \monprod N \ar@<-.5ex>[r] \ar@<.5ex>[r]
&M \monprod N \ar[r] &\strut M \balprod N, }
\]
where one of the lefthand arrows is induced by the right action of
$\monoid$ on $M$ and the other by the left action of $\monoid$ on $N$.


It is clear from the definitions that the set $\Leftmod(L,M)$ of left
$\monoid$-module maps from $L$ to $M$ is naturally in bijective
correspondence with the set $\Moncat(\unit,\lmfunc{L}{M})$ of maps in
$\Moncat$ from the unit $\unit$ to $\lmfunc{L}{M}$.  Thus, the left
$\monoid$-function object construction enriches the category
$\Leftmod$ of left $\monoid$-modules over the category \Moncat.  For
objects $X$ in $\Moncat$ and $L,M$ in $\Leftmod$, the objects
$L\monprod X$ and $\func{X}{M}$ inherit left $\monoid$-actions from
$L$ and the enriched parametrized adjunctions
\[
\lmfunc{L\monprod X}{M}\iso \func{X}{\lmfunc{L}{M}}\iso
\lmfunc{L}{\func{X}{M}}
\]
indicate that the constructions $L\monprod X$ and $\func{X}{M}$
provide tensors and cotensors for the enrichment of $\Leftmod$ over
\Moncat.  The forgetful functor \vmap{\Leftmod}{\Moncat} has enriched
left and right adjoints, called the free and cofree functors, sending
an object $X$ of $\Moncat$ to $A \monprod X$ and $\func{A}{X}$,
respectively.  These functors have a rich structure of interrelations
and coherences that the enriched category theory language concisely
encodes and which would be tedious to list in terms of individual
natural isomorphisms.

Our first objective is to describe conditions under which all of this
structure passes over to the homotopy categories.  Much of it passes
over with no restrictions other than a very standard one on the model
structure inherited by \Leftmod from \Moncat.  A closed model
structure on the module category \Leftmod is said to have fibrations
and weak equivalences created in \Moncat if a map $f$ in \Leftmod is a
fibration or weak equivalence in the model structure for \Leftmod if
and only if it is one in the model structure for \Moncat.  The
following is the most basic theorem in this direction.

\begin{theorem}\label{intenr}
Let $\Moncat$ be a monoidal model category, and let \monoid be a
monoid in $\Moncat$.  If the category of left $\monoid$-modules is a
closed model category with fibrations and weak equivalences created in
$\Moncat$, then:
\begin{thmlist}

\item The right derived functor $\emptyderlmfunc$ of $\lmfunc{-}{-}$
exists and enriches $\Ho\Leftmod$ over $\Ho\Moncat$.

\item The right derived functor $\emptyderfunc$ of the cotensor
functor $\func{-}{-}\colon \Moncat^{\op}\times \Leftmod\to \Leftmod$
exists and provides cotensors for $\Ho\Leftmod$ over $\Ho\Moncat$

\item The left derived functor $\emptydermonprod$ of the tensor
functor $\monprod\colon\Leftmod\times \Moncat \to \Leftmod$ exists and
provides tensors for $\Ho\Leftmod$ over $\Ho\Moncat$
\end{thmlist}
\end{theorem}

The enrichment concisely encodes many relations and coherences that
are less obvious for these derived functors than for the corresponding
functors on $\Moncat$ and $\Leftmod$.  For example, the interpretation
of $\emptydermonprod$ as a tensor encodes coherent associativity
natural isomorphisms as well as various adjunctions.  This theorem is
a special case of a general theorem for closed model categories
enriched over monoidal model categories, discussed in
Section~\ref{EnrModCatSec}.

The condition that fibrations and weak equivalences are created in
$\Moncat$ obviously implies that the forgetful functor from \Leftmod
to \Moncat and its left adjoint free functor form a Quillen
adjunction, and so induce a derived adjunction on the homotopy
categories.  Likewise, since the cotensor in $\Ho\Leftmod$ is the
right derived functor $\emptyderfunc$ of $\func{-}{-}\colon
\Moncat^{\op}\times \Leftmod\to \Leftmod$, it follows that its
composition with the derived forgetful functor to $\Ho\Moncat$ is
naturally isomorphic to the right derived functor of
$\func{-}{-}\colon \Moncat^{\op}\times \Moncat\to \Moncat$.  The
corresponding assertions about the existence of a right adjoint for
the derived forgetful functor and that the derived forgetful functor
preserves tensors need not hold in general, but require an additional
hypothesis on $\monoid$.  This hypothesis depends only on $\monoid$
viewed as an object of $\Moncat$, and it is convenient to state in a
general context for further use in the statements below.

\begin{definition}
Let $\EC$ be a closed model category that is also enriched over a
monoidal model category $\Moncat$ by function objects $\efunc{-}{-}$.
An object $C$ is said to be \defterm{semicofibrant} in $\EC$ when the functor
$\efunc{C}{-}$ preserves fibrations and acyclic fibrations.
\end{definition}

This property is explored further in Section~\ref{SecSemi}.  The
following result proved in that section provides enough information
about this notion for our present purposes.

\begin{proposition}\label{SemiCoPropA}
Let \monoid be a monoid in a monoidal model category $\Moncat$ for
which the module category \Leftmod is a closed model category with
fibrations and weak equivalences created in $\Moncat$.
\begin{thmlist} 

\item $M$ is semicofibrant in \Leftmod if and only if the
functor \map{M\monprod(-)}{\Moncat}{\Leftmod} preserves cofibrations
and acyclic cofibrations.

\item If $M \in \Leftmod$ is cofibrant in \Leftmod then it is
semicofibrant in \Leftmod. Moreover, if the unit \unit is cofibrant in
\Moncat, then an object $M$ of \Leftmod is semicofibrant in \Leftmod
if and only if it is cofibrant in \Leftmod.

\item\llabel{intsemicof} If $M \in \Leftmod$ is semicofibrant in
\Leftmod and $M\to N$ is a cofibration in \Leftmod, then $N$ is
semicofibrant in \Leftmod.

\item\llabel{intsemiunit} \monoid, considered as an object of
\Leftmod, is semicofibrant in \Leftmod.  In particular, $\unit$ is
semicofibrant in \Moncat.

\end{thmlist}
\end{proposition}

If the monoid \monoid is semicofibrant when considered as an object of
\Moncat, then all of the enriched structure of \Leftmod over \Moncat
discussed above passes to the homotopy categories.

\begin{theorem}\label{ForgetAdjProp}
Let $\Moncat$ be a monoidal model category, and let \monoid be a
monoid in $\Moncat$.  If the category of left $\monoid$-modules is a
closed model category with fibrations and weak equivalences created in
$\Moncat$, then:
\begin{thmlist}

\item The left derived functor of the free functor
\vmap{\Moncat}{\Leftmod} exists and is a $\Ho\Moncat$-enriched left
adjoint to the forgetful functor from $\Ho\Leftmod$ to $\Ho\Moncat$.

\item The forgetful functor \vmap{\Ho\Leftmod}{\Ho\Moncat} preserves
cotensors.  \par\smallskip\noindent Moreover, if $\monoid$ is
semicofibrant as an object of \Moncat, then: \smallskip

\item The right derived functor of the cofree functor
\vmap{\Moncat}{\Leftmod} exists and is a $\Ho\Moncat$-enriched right
adjoint to the forgetful functor from $\Ho\Leftmod$ to $\Ho\Moncat$.

\item The forgetful functor \vmap{\Ho\Leftmod}{\Ho\Moncat} preserves
tensors.
\end{thmlist}
\end{theorem}

The hypothesis above that $\monoid$ is semicofibrant as an object of
$\Moncat$ seems to hold quite generally: Often the category of monoids
in $\Moncat$ forms a closed model category where the unit map $I\to
\monoid$ for a cofibrant monoid $\monoid$ is a cofibration in
$\Moncat$.  In that case parts~(\ref{intsemicof})
and~(\ref{intsemiunit}) of Proposition~\ref{SemiCoPropA} imply that
cofibrant monoids are semicofibrant objects of $\Moncat$.  This
applies in particular when the hypotheses of the main theorem of
Schwede--Shipley \cite[4.1]{ssmonoidal} hold.

The results above are special cases of more general results about
bimodule categories and functors between such categories.  If
\monoidone and \monoidtwo are monoids in $\Moncat$, then an
$(\monoidone,\monoidtwo)$-bimodule $M$ in \Moncat is an object of
\Moncat with commuting left \monoidone-module and right
\monoidtwo-module structures.  Equivalently, it may be described as a
left $\monoidone \monprod \monoidtwo^\op$-module.  The category of
$(\monoidone,\monoidtwo)$-bimodules is denoted
\Bimod{\monoidone}{\monoidtwo}, and the $\monoidone \monprod
\monoidtwo^\op$-module function object $\lmfunc[\monoidone \monprod
\monoidtwo^\op]{-}{-}$ is denoted
$\bmfunc{\monoidone}{\monoidtwo}{-}{-}$.

Any object $X$ of \Moncat carries a canonical $(\unit,\unit)$-bimodule
structure.  Moreover, any right \monoidone-module $N$ carries a
canonical $(\unit,\monoidone)$-bimodule structure.  Analogously, any
left \monoidtwo-module $M$ is canonically a
$(\monoidtwo,\unit)$-bimodule.  Thus, we can identify the category
\Moncat with the categories \Leftmod[\unit], \Rightmod[\unit] and
\Bimod{\unit}{\unit}.  Similarly, \Leftmod[\monoidone] and
\Rightmod[\monoidtwo] can be identified with \Bimod{\monoidone}{\unit}
and \Bimod{\unit}{\monoidtwo}, respectively.  With this perspective
all results below for categories of bimodules specialize to
corresponding results for categories of left and/or right modules.

For any $(\monoidone,\monoidtwo)$-bimodule $M$ and any left
\monoidtwo-module $N$, the balanced product $M \monprod_{\monoidtwo}
N$ is naturally a left \monoidone-module.  More generally, for monoids
\monoidone, \monoidtwo, and \monoidthree, we can consider
$-\monprod_{\monoidtwo}-$ to be a functor
\begin{equation*}
\map{-\monprod_{\monoidtwo}-}{\Bimod{\monoidone}{\monoidtwo} \times
\Bimod{\monoidtwo}{\monoidthree}}{\Bimod{\monoidone}{\monoidthree}}.
\end{equation*}
Similarly, if $M$ is an $(\monoidone,\monoidtwo)$-bimodule and $P$ is a
right \monoidone-module, then the function objects $\lmfunc{M}{P}$ and
$\lmfunc{P}{M}$ inherit \monoidtwo-actions from $M$.  These actions
are left and right, respectively.  Generalizing this, we can think of
the left \monoidone-module function object construction
$\lmfunc{-}{-}$ as a functor
\begin{equation*}
\lmfunc[\monoidone]{-}{-}\colon
\Bimod{\monoidone}{\monoidtwo}^{\op}\times
\Bimod{\monoidone}{\monoidthree}\to \Bimod{\monoidtwo}{\monoidthree}.
\end{equation*}
Analogously, we can think of the right \monoidthree-module function
object construction $\rmfunc[\monoidthree]{-}{-}$ as a functor
\begin{equation*}
\rmfunc[\monoidthree]{-}{-}\colon
\Bimod{\monoidone}{\monoidthree}^{\op}\times
\Bimod{\monoidtwo}{\monoidthree}\to \Bimod{\monoidtwo}{\monoidone}.
\end{equation*}
The three constructions $\lmfunc[\monoidone]{-}{-}$,
$-\monprod_{\monoidtwo}-$, and $\rmfunc[\monoidthree]{-}{-}$ are
related by natural isomorphisms
\begin{equation}\label{BimHomTenAdj}
\bmfunc{\monoidone}{\monoidtwo}{M}{\rmfunc[\monoidthree]{N}{P}} \iso
\bmfunc{\monoidone}{\monoidthree}{M\balprod[\monoidtwo]N}{P} \iso
\bmfunc{\monoidtwo}{\monoidthree}{N}{\lmfunc[\monoidone]{M}{P}},
\end{equation}
for $M$ in $\Bimod{\monoidone}{\monoidtwo}$, $N$ in
$\Bimod{\monoidtwo}{\monoidthree}$, and $P$ in
$\Bimod{\monoidone}{\monoidthree}$.  Note that the second isomorphism
is precisely the first isomorphism for the opposite monoids under the
isomorphisms of categories
\begin{equation}\label{catiso}
\Bimod{\monoidtwo}{\monoidthree}\iso
\Bimod{\monoidthree^{\op}}{\monoidtwo^{\op}}\qquad
\Bimod{\monoidone}{\monoidtwo}\iso
\Bimod{\monoidtwo^{\op}}{\monoidone^{\op}}\qquad
\Bimod{\monoidone}{\monoidthree}\iso
\Bimod{\monoidthree^{\op}}{\monoidone^{\op}}.
\end{equation}

All of our module categories are enriched over \Moncat and each of the
constructions $\lmfunc[\monoidone]{-}{-}$, $-\monprod_{\monoidtwo}-$,
and $\rmfunc[\monoidthree]{-}{-}$ gives a \Moncat-enriched bifunctor.
With this viewpoint, the two isomorphisms above become two enriched
parametrized adjunctions.  The following theorem extends these
adjunctions to the homotopy categories.

\begin{theorem}\label{intexist}
Let $\Moncat$ be a monoidal model category, and let $\monoidone$,
$\monoidtwo$, and $\monoidthree$ be monoids in $\Moncat$.  Assume that
each category of modules in each statement below is a closed model
category with fibrations and weak equivalences created in $\Moncat$.
If $\monoidtwo$ is semicofibrant when considered as an object in
$\Moncat$, then:
\begin{thmlist}
\item The right derived functor
$\WowExt{\monoidone}{\monoidtwo}{\monoidthree}{-}{-}$ of
\[
\lmfunc[\monoidone]{-}{-}\colon
\Bimod{\monoidone}{\monoidtwo}^{\op}\times
\Bimod{\monoidone}{\monoidthree}\to \Bimod{\monoidtwo}{\monoidthree}
\]
exists and is enriched over $\Ho\Moncat$.
\item\label{IEintwowtor}The left derived functor
$\WowTor{\monoidone}{\monoidtwo}{\monoidthree}{-}{-}$ of
\[
\balprod[\monoidtwo]\colon \Bimod{\monoidone}{\monoidtwo}\times
\Bimod{\monoidtwo}{\monoidthree}\to \Bimod{\monoidone}{\monoidthree}
\]
exists, is enriched over $\Ho\Moncat$, and forms an enriched
parametrized adjunction with
$\WowExt{\monoidone}{\monoidtwo}{\monoidthree}{-}{-}$:
\begin{equation}\label{HoBimHomTenAdjA}
\derbmfunc{\monoidone}{\monoidthree}{\WowTor{\monoidone}{\monoidtwo}{\monoidthree}{M}{N}}{P}\iso
\derbmfunc{\monoidtwo}{\monoidthree}{N}{\WowExt{\monoidone}{\monoidtwo}{\monoidthree}{M}{P}}.
\end{equation}
\item\llabel{intquiladj} Let $M$ be a
$(\monoidone,\monoidtwo)$-bimodule.  If the underlying left
$\monoidone$-module of $M$ is semicofibrant in \Leftmod, then
$\WowTor{\monoidone}{\monoidtwo}{\monoidthree}{M}{-}$ is the left
derived functor of $M\balprod[\monoidtwo]{-}$,
$\WowExt{\monoidone}{\monoidtwo}{\monoidthree}{M}{-}$ is the right
derived functor of $\lmfunc[\monoidone]{M}{-}$, and the adjunction
\begin{equation*}
\xymatrix{{M\balprod[\monoidtwo](-):\
\Bimod{\monoidtwo}{\monoidthree}\ar@<0.7ex>[r]} &
{\Bimod{\monoidone}{\monoidthree}\ :
\lmfunc[\monoidone]{M}{-}\ar@<0.7ex>[l]}}
\end{equation*}
is a Quillen adjunction.
\end{thmlist}
\end{theorem}  

See Definition~\ref{totalderivedfunctor} in Section~\ref{BiFunSec} for
a precise definition of the enrichment of the derived functor of a
bifunctor.

Part~(\ref{intquiladj}) above applies in particular when $M$ is a
cofibrant $(\monoidone,\monoidtwo)$-bimodule, because then the
underlying $\monoidone$-module of $M$ is cofibrant in $\Leftmod$ and
therefore semicofibrant in $\Leftmod$.  To see this, note that the
right adjoint of the forgetful functor
\vmap{\Bimod{\monoidone}{\monoidtwo}}{\Leftmod} is the functor
$\func{\monoidtwo}{-}\colon \Leftmod\to
\Bimod{\monoidone}{\monoidtwo}$, which preserves fibrations and
acyclic fibrations since by hypothesis $\monoidtwo$ is semicofibrant
in $\Moncat$.  It follows that the forgetful functor
\vmap{\Bimod{\monoidone}{\monoidtwo}}{\Leftmod} preserves cofibrations
and acyclic cofibrations.

Another interesting case of part~(\ref{intquiladj}) occurs when
$\monoidone=\monoidtwo$ and $M=\monoidtwo$.  Then the statement is
that the functors
$\WowTor{\monoidtwo}{\monoidtwo}{\monoidthree}{\monoidtwo}{-}$ and
$\WowExt{\monoidtwo}{\monoidtwo}{\monoidthree}{\monoidtwo}{-}$ are
naturally isomorphic to the identity functor on
$\Ho\Bimod{\monoidtwo}{\monoidthree}$.

Isomorphism \eqref{HoBimHomTenAdjA} of Theorem \ref{intexist} can be
coupled with the isomorphisms \eqref{catiso} of categories to obtain
the pair of enriched isomorphisms
\begin{equation}\label{HoBimHomTenAdjB}
\mbox{\small $%
\derbmfunc{\monoidone}{\monoidtwo}{M}{\WowExtP{\monoidone}{\monoidtwo}{\monoidthree}{N}{P}}
\iso
\derbmfunc{\monoidone}{\monoidthree}{\WowTor{\monoidone}{\monoidtwo}{\monoidthree}{M}{N}}{P}\iso\derbmfunc{\monoidtwo}{\monoidthree}{N}{\WowExt{\monoidone}{\monoidtwo}{\monoidthree}{M}{P}}
$}
\end{equation} 
which are the derived versions of the isomorphisms
\eqref{BimHomTenAdj}.

The universal property of derived functors implies that functors in
the previous theorem are appropriately natural in the monoids
$\monoidone$, $\monoidtwo$, $\monoidthree$.  In fact, the natural
transformations so obtained are enriched.

\begin{theorem}\label{intnat}
Let $\Moncat$ be a monoidal model category, and assume that each
category of modules in each statement below is a closed model category
with fibrations and weak equivalences created in $\Moncat$.
\begin{thmlist}
\item When the underlying objects of $\monoidtwo$ and $\monoidtwo'$
are semicofibrant in $\Moncat$, maps of monoids $\monoidone'\to
\monoidone$, $\monoidtwo'\to \monoidtwo$, and $\monoidthree'\to
\monoidthree$, induce an enriched natural transformation of bifunctors
(from
$\Ho\Bimod{\monoidone}{\monoidtwo}^{\op},\Ho\Bimod{\monoidone}{\monoidthree}$
to $\Ho\Bimod{\monoidtwo'}{\monoidthree'}$)
\[
\WowExt{\monoidone}{\monoidtwo}{\monoidthree}{-}{-}
\to\WowExt{\monoidone'}{\monoidtwo'}{\monoidthree'}{-}{-}
\]
making $\WExt$ appropriately functorial in the monoid variables.  In
particular, this transformation is compatible with the natural
transformation
$\lmfunc[\monoidone]{-}{-}\to\lmfunc[\monoidone']{-}{-}$,
\item When the underlying objects of $\monoidtwo$ and $\monoidtwo'$
are semicofibrant in $\Moncat$, maps of monoids $\monoidone'\to
\monoidone$, $\monoidtwo'\to \monoidtwo$, and $\monoidthree'\to
\monoidthree$, induce an enriched natural transformation of bifunctors
(from
$\Ho\Bimod{\monoidone}{\monoidtwo},\Ho\Bimod{\monoidtwo}{\monoidthree}$
to $\Ho\Bimod{\monoidone'}{\monoidthree'}$)
\[
\WowTor{\monoidone'}{\monoidtwo'}{\monoidthree'}{-}{-}
\to\WowTor{\monoidone}{\monoidtwo}{\monoidthree}{-}{-},
\]
making $\WTor$ appropriately functorial in the monoid variables.  In
particular, this transformation is compatible with the natural
transformation $\balprod[\monoidtwo']\to \balprod[\monoidtwo]$,

\item The $\WExt$ and $\WTor$ natural transformations above are
appropriately compatible with the adjunction isomorphism of
Theorem~\ref{IEintwowtor}.
\end{thmlist}
\end{theorem}

In favorable situations, the underlying object in $\Ho\Moncat$ of
$\WowTor{\monoidone}{\monoidtwo}{\monoidthree}{-}{-}$ should only
depend on $\monoidtwo$ and not on $\monoidone$ and $\monoidthree$.
Similarly, the underlying object in $\Ho\Moncat$ of
$\WowExt{\monoidone}{\monoidtwo}{\monoidthree}{-}{-}$ should only
depend on $\monoidone$ and not on $\monoidtwo$ and $\monoidthree$.
The natural transformations of Theorem \ref{intnat} allow us to
convert this intuition into the following precise statement.  In it,
we drop the notation for any monoid variable when it is the unit
$\unit$.

\begin{theorem}\label{intforcomp}
Let $\Moncat$ be a monoidal model category, and let $\monoidone$,
$\monoidtwo$, and $\monoidthree$ be monoids in $\Moncat$.  Assume that
each category of modules in each statement below is a closed model
category with fibrations and weak equivalences created in $\Moncat$.
If $\monoidtwo$ is semicofibrant when considered as an object in
$\Moncat$, then:

\begin{thmlist}
\item\llabel{intforgext} The natural transformation
\[
\WowExt{\monoidone}{\monoidtwo}{\monoidthree}{-}{-}\to
\WowExt{\monoidone}{}{}{-}{-}= \derlmfunc[\monoidone]{-}{-}
\]
in $\Ho\Moncat$ induced by the unit maps $\unit\to\monoidtwo$ and
$\unit\to \monoidthree$ is an isomorphism.
\item\llabel{intforgtor} If the underlying object of $\monoidone$ is
semicofibrant in $\Moncat$, then the natural transformation
\[
\WowTor{}{\monoidtwo}{\monoidthree}{-}{-} \to
\WowTor{\monoidone}{\monoidtwo}{\monoidthree}{-}{-}
\]
in $\Ho\Rightmod[\monoidthree]$ induced by the unit map
$\unit\to\monoidone$ is an isomorphism.  Similarly, if the underlying
object of $\monoidthree$ is semicofibrant in $\Moncat$, then the
natural transformation
\[
\WowTor{\monoidone}{\monoidtwo}{}{-}{-} \to
\WowTor{\monoidone}{\monoidtwo}{\monoidthree}{-}{-}
\]
in $\Ho\Rightmod[\monoidone]$ induced by the unit map
$\unit\to\monoidthree$ is an isomorphism.

\end{thmlist}
\end{theorem}

Since a map in $\Ho\Bimod{\monoidtwo'}{\monoidthree'}$ is an
isomorphism if and only if it is sent to an isomorphism in
$\Ho\Moncat$, it follows from part~(\ref{intforgext}) that for any
maps of monoids $\monoidtwo'\to \monoidtwo$ and $\monoidthree\to
\monoidthree'$, when $\monoidtwo$ and $\monoidtwo'$ have underlying
objects that are semicofibrant in $\Moncat$, the induced map
$\WowExt{\monoidone}{\monoidtwo}{\monoidthree}{-}{-}\to
\WowExt{\monoidone}{\monoidtwo'}{\monoidthree'}{-}{-}$ is an
isomorphism in $\Ho\Bimod{\monoidtwo'}{\monoidthree'}$. Likewise, it
follows from part~(\ref{intforgtor}) that for any map of monoids
$\monoidone'\to\monoidone$ whose underlying objects are semicofibrant
in $\Moncat$, the induced map
$\WowTor{\monoidone'}{\monoidtwo}{\monoidthree}{-}{-}\to
\WowTor{\monoidone}{\monoidtwo}{\monoidthree}{-}{-}$ is an isomorphism
in $\Ho\Bimod{\monoidone}{\monoidthree}$.  Similarly, for any map of
monoids $\monoidthree'\to\monoidthree$ whose underlying objects are
semicofibrant in $\Moncat$, the induced map
$\WowTor{\monoidone}{\monoidtwo}{\monoidthree'}{-}{-}\to
\WowTor{\monoidone}{\monoidtwo}{\monoidthree}{-}{-}$ is an isomorphism
in $\Ho\Bimod{\monoidone}{\monoidthree}$.

Since for Quillen adjunctions the left derived functor of the
composite of the left adjoints is the composite of the left derived
functors, the last part of Theorem \ref{intexist} gives an
``associativity'' isomorphism for the derived functors.

\begin{theorem}\label{intassoc}
Let $\Moncat$ be a monoidal model category, and let $\monoidone$,
$\monoidtwo$, $\monoidthree$, and $\monoidfour$ be monoids in
$\Moncat$.  Assume that each category of modules in the statement
below is a closed model category with fibrations and weak equivalences
created in $\Moncat$.  If the underlying objects of $\monoidtwo$ and
$\monoidthree$ are semicofibrant in $\Moncat$, then there is a
canonical enriched natural isomorphism of trifunctors
\[
\WowTor{\monoidone}{\monoidtwo}{\monoidfour}{L}{%
(\WowTor{\monoidtwo}{\monoidthree}{\monoidfour}{M}{N})}\iso
\WowTor{\monoidone}{\monoidthree}{\monoidfour}{%
(\WowTor{\monoidone}{\monoidtwo}{\monoidthree}{L}{M})}{N},
\]
compatible with the associativity isomorphism for the symmetric
monoidal product in $\Ho\Moncat$ and satisfying the evident analogue
of the pentagon law.  Adjointly, there is a canonical enriched natural
isomorphism of trifunctors
\[
\WowExt{\monoidone}{\monoidthree}{\monoidfour}{%
(\WowTor{\monoidone}{\monoidtwo}{\monoidthree}{M}{N})}{P}\iso
\WowExt{\monoidtwo}{\monoidthree}{\monoidfour}{N}{%
(\WowExt{\monoidone}{\monoidtwo}{\monoidfour}{M}{P})}.
\]
\end{theorem}

Several special cases of the results presented above are of particular
interest.  These include:

\subsection*{Tensors and cotensors}
Although treated explicitly in Theorem~\ref{intenr}, the existence and
interpretation of tensors and cotensors as derived functors also
follows from the general bimodule theorems above.  Tensors and
cotensors in \Leftmod comprise the special case of the isomorphisms
\eqref{BimHomTenAdj} in which $\monoidtwo = \monoidthree = \unit$.
Likewise, tensors and cotensors in $\Ho\Leftmod$ comprise the special
case of the isomorphisms \eqref{HoBimHomTenAdjB} in which $\monoidtwo
= \monoidthree = \unit$.  This indicates that
$\WowTor{\monoidone}{}{}{M}{X}$ provides the tensor $\dertensor{M}{X}$
for $\Ho\Leftmod[\monoidone]$.  The last part of the
Theorem~\ref{intforcomp} indicates that the tensor $\dertensor{M}{X}$
in $\Ho\Leftmod[\monoidone]$ agrees with the derived monoidal product
$\dermonprod{M}{X}$ in $\Ho\Moncat$ when the underlying object in
$\Moncat$ of $\monoidone$ is semicofibrant.  Moreover, it follows that
for a map of monoids $\monoidone\to \monoidtwo$ whose underlying
objects are semicofibrant, the derived forgetful (or ``pullback'')
functor $\Ho\Leftmod[\monoidtwo]\to \Ho\Leftmod[\monoidone]$ preserves
tensors.  This special case of isomorphisms \eqref{HoBimHomTenAdjB}
also implies that $\WowExtp{\monoidone}{}{}{X}{M}$ provides the
cotensors for $\Ho\Leftmod$, and that these are preserved by the
derived forgetful functor to $\Ho\Moncat$.  Note that tensors are
preserved by all enriched left adjoints and cotensors are preserved by
all enriched right adjoints, and so the remarks on preservation of
tensors and cotensors also follow from the observations on extension
of scalars and coextension of scalars below.

\subsection*{Extension of scalars}
Let $\monoidtwo\to \monoidone$ be a map of monoids in $\Moncat$, and
assume that the categories $\Leftmod[\monoidone]$ and
$\Leftmod[\monoidtwo]$ are closed model categories with fibrations and
weak equivalences created in $\Moncat$.  It then follows formally that
the extension of scalars functor
$\monoidone\balprod[\monoidtwo](-)\colon \Leftmod[\monoidtwo]\to
\Leftmod[\monoidone]$ and the forgetful (or pullback) functor
$\Leftmod[\monoidone]\to \Leftmod[\monoidtwo]$ form a Quillen
adjunction.  When $\monoidone$ and $\monoidtwo$ are semicofibrant in
$\Moncat$, Theorem~\ref{intexist} implies that the left derived
extension of scalars functor is given by
$\WowTor{\monoidone}{\monoidtwo}{}{\monoidone}{-}$ and
Theorem~\ref{intforcomp} implies that it is naturally isomorphic to
$\WowTor{}{\monoidtwo}{}{\monoidone}{-}$ in $\Ho\Moncat$.  In
particular, in this case, when $\monoidtwo \to \monoidone$ is a weak
equivalence, the extension of scalars adjunction is a Quillen
equivalence.

\subsection*{Coextension of scalars} Let $\monoidone\to \monoidtwo$ be
a map of monoids in $\Moncat$, and assume that the categories
$\Leftmod[\monoidone]$ and $\Leftmod[\monoidtwo]$ are closed model
categories with fibrations and weak equivalences created in $\Moncat$.
The forgetful functor $\Leftmod[\monoidtwo]\to \Leftmod[\monoidone]$
has a right adjoint given by $\lmfunc[\monoidone]{\monoidtwo}{-}$.  If
$\monoidtwo$ is semicofibrant in $\Leftmod[\monoidone]$, then
Theorem~\ref{intexist} implies that this is a Quillen adjunction and
it identifies the right derived coextension of scalars functor as
$\WowExt{\monoidone}{\monoidtwo}{}{\monoidtwo}{-}$.  Moreover,
Theorem~\ref{intforcomp} implies that this functor is naturally
isomorphic to $\derlmfunc[\monoidone]{\monoidtwo}{-}$ in $\Ho\Moncat$.

\subsection*{Free and cofree functors bimodule structure}
For the map of monoids $\unit\to \monoid$, the extension of scalars
functor and coextension of scalars functor are called the free functor
and the cofree functor.  These functors have the extra structure that
they factor through the forgetful functor $\Bimod{\monoid}{\monoid}\to
\Leftmod[\monoid]$.  Assume that the categories of $\Leftmod$ and
$\Bimod{\monoid}{\monoid}$ are closed model categories with fibrations
and weak equivalences created in $\Moncat$ and that $\monoid$ is
semicofibrant in $\Moncat$.  Then $\monoid\monprod \monoid^{\op}$ is
also semicofibrant in $\Moncat$ and Theorem~\ref{intforcomp} implies
that the functors
\[
\derfree (-)=\WowTor{\monoid\monprod\monoid^{\op}}{}{}{\monoid}{-},
\tand
\dercofree(-)=\WowExt{}{\monoid\monprod\monoid^{\op}}{}{\monoid}{-},
\]
provide factorizations of the derived free and cofree functors through
$\Ho\Bimod{\monoid}{\monoid}$.  The adjunctions also identify
$\derfree X$ as $\dertensor{A}{X}$, the tensor (in
$\Ho\Bimod{\monoid}{\monoid}$) of $\monoid$ with the object $X$ of
$\Ho\Moncat$.  Since tensors commute with enriched left adjoints, we
obtain natural isomorphisms in $\Ho\Leftmod$ (or
$(\Ho\Leftmod)^{\op}$)
\[
\dertensor{M}{X}\iso \WowTor{\monoid}{\monoid}{}{\derfree X}{M}, \tand
\dercotensor{X}{M}\iso \WowExt{\monoid}{\monoid}{}{\derfree X}{M}.
\]
The isomorphism $\derlmfunc{M}{\dercofree X}\iso \derfunc{M}{X}$ from
the enriched adjunction also refines to an isomorphism in
$\Ho\Rightmod$
\[
\WowExt{\monoid}{}{\monoid}{M}{\dercofree X}\iso
\WowExtR{}{\monoid}{}{M}{X}
\]
as an instance of the universal map of enriched derived functors.
Although not a direct result of the results listed above, this last
enriched natural transformation is an immediate consequence of the
more general Theorem~\ref{penrquilcri} in Section~\ref{BiFunSec}
below.

\medskip

In practice, many monoidal model categories have additional properties
that make the semicofibrant hypotheses in the results above
unnecessary in certain cases.  The process of eliminating these
hypotheses is discussed in Section~\ref{NonSemiSec}.

\medskip

The second author would like to thank Brooke Shipley and Andrew
Blumberg for helpful comments. 


\section{Monoidal model categories}\label{MndlModCatSec}

This section reviews the terminology and basic theory of monoidal
model categories from \cite{ssmonoidal} (and \cite{HvModCat}).  The
definition of a monoidal model category involves constraints on the
interaction of the model structure with the closed symmetric monoidal
structure.  The imposed conditions suffice to ensure that the homotopy
category inherits a closed symmetric monoidal category structure and
that the localization functor is lax symmetric monoidal.  The
conditions are stated in terms of the following two standard maps.
Let \map{f}{A}{B}, \map{g}{K}{L}, and \map{h}{X}{Y} be maps in a
symmetric monoidal closed category $\Moncat$.  Then the maps
\[
\boxfunc{g}{h}\colon \func{L}{X}\to
\func{K}{X}\times_{\func{K}{Y}}\func{L}{Y}
\]
and
\[
f\boxmonprod g\colon (A\monprod L) \cup_{A\monprod K}(B\monprod K)\to
B\monprod L
\]
are defined by the diagrams
\begin{equation}\label{diaboxfunc}
\vcenter{
\xymatrix@-1ex{
{\func{L}{X} \ar[drr]^-{\boxfunc{g}{h}} \ar@/^1pc/[drrr] \ar@/_1pc/[ddrr]} & {} & {} & {} \\
{} & {} & {\func{K}{X}\times_{\func{K}{Y}}\func{L}{Y} \ar[r] \ar[d]}
	& {\func{L}{Y} \ar[d]^-{g^*}}  \\
{} & {} & {\func{K}{X} \ar[r]_-{h_*}} & {\func{K}{Y}} }}
\end{equation}
and
\begin{equation}\label{diaboxprod}
\vcenter{
\xymatrix@-1ex{ {A \monprod K \ar[r]^-{f \monprod \, \id} \ar[d]_-{\id
\monprod \, g}}
	& {B \monprod  K \ar[d] \ar@/^1pc/[ddrr]} & {} & {}\\
{A \monprod L \ar[r] \ar@/_1pc/[drrr]} & {(A\monprod L)
\cup_{A\monprod K}(B\monprod K) \ar[drr]^-{f\,\boxmonprod\, g}}
	& {} & {}\\
{} & {} & {} & {B \monprod L} }}
\end{equation}
in which the squares are a pullback and a pushout, respectively.

\begin{definition}\label{defmmc}
A \defterm{monoidal model category} $\Moncat$ is a closed model category with a
closed symmetric monoidal structure satisfying the following two
axioms:
\begin{description}
\item[\ENa] If $g\colon K\to L$ is a cofibration and $h\colon X\to Y$
is a fibration, then $\boxfunc{g}{h}$ is a fibration.  Moreover, if
either $g$ or $h$ is also a weak equivalence, then so is
$\boxfunc{g}{h}$.

\item[\UF] There exists a cofibrant object $\cofunit$ and a weak
equivalence $\omega \colon \cofunit \to \unit$ such that the composite
$\cuf$ of the adjoint of the unit isomorphism and $\omega^{*}$
\[
\cuf\colon Z\to \func{\unit}{Z}\to \func{\cofunit}{Z}
\]
is a weak equivalence for every fibrant object $Z$.
\end{description}
\end{definition}

The first axiom, the \defterm{Enrichment Axiom} is the internal version of
Quillen's axiom \SMa.  We have given it in a form that easily
generalizes to the context of enriched categories in the next section.
Each of the above axioms may be reformulated adjointly in terms of
$\monprod$, and these reformulations seem to be easier to work with in
practice.  The adjoint form of the Enrichment Axiom is called the
\defterm{Pushout Product Axiom} \cite[3.1]{ssmonoidal}

\begin{proposition}\label{enrpopr}
Let $\Moncat$ be a closed model category with a closed symmetric
monoidal structure.  Then $\Moncat$ satisfies the Enrichment Axiom
\ENa if and only if it satisfies the Pushout Product Axiom:
\begin{description}
\item[\PPa] If $f\colon A\to B$ and $g\colon K\to L$ are cofibrations,
then so is $f\boxmonprod g$. Moreover, if either $f$ or $g$ is also a
weak equivalence, then so is $f\boxmonprod g$.
\end{description}
When $\Moncat$ satisfies these axioms, it satisfies the Unit Axiom \UF
if and only if it satisfies the following axiom:
\begin{description}
\item[\UC] There exists a cofibrant object $\cofunit$ and a weak
equivalence $\omega \colon \cofunit \to \unit$ such that the composite
$\cum$ of $\id\monprod\omega$ and the unit isomorphism
\[
\cum\colon X\monprod \cofunit \to X\monprod \unit\iso X
\]
is a weak equivalence for every cofibrant object $X$.
\end{description}
\end{proposition}

The equivalence of the axioms \ENa and \PPa follows from the
characterization of (acyclic) cofibrations and (acyclic) fibrations in
terms of lifting properties, using the
$(-\monprod-,\emptyfunc)$-adjunction applied to $f$, $g$, and $h$ as
above (see, for example, \cite[4.2.2]{HvModCat}).  The equivalence of
the two unit axioms is closely related to the construction of the
derived product and function functors and our discussion of it is
postponed until after Proposition~\ref{derivedbifunctors} below.

To describe the implications of \PPa for the functors $\monprod $ and
$\emptyfunc$, we must first recall the standard model structures on
the opposite of a closed model category and the product of two closed
model categories.  If $\Moncat$ and $\Moncat'$ are closed model
categories, then $\Moncat' \times \Moncat$ is a closed model category
whose cofibrations, fibrations, and weak equivalences are the maps
that are cofibrations, fibrations, and weak equivalences,
respectively, in each coordinate.  Also, $\Moncat^{\op}$ is a closed
model category whose cofibrations, fibrations, and weak equivalences
are the maps that are the opposites of fibrations, cofibrations, and
weak equivalences, respectively.  In particular, the fibrant objects
in $\Moncat^{\op}$ are the cofibrant objects in $\Moncat$.

\begin{proposition}\label{pppweq}
Let $\Moncat$ be a monoidal model category, and $X$, $Z$ objects of
$\Moncat$ that are cofibrant and fibrant, respectively.
\begin{thmlist}
\item The functors $X \monprod (-)$ and $(-)\monprod X$ preserve
cofibrations and acyclic cofibrations.

\item The functor $\func{X}{-}$ preserves fibrations and acyclic
fibrations.

\item The functor $\func{-}{Z}$ converts cofibrations and acyclic
cofibrations in $\Moncat$ into fibrations and acyclic fibrations,
respectively.

\item The functor $\monprod \colon \Moncat\times \Moncat\to \Moncat$
preserves weak equivalences between cofibrant objects

\item The functor $\emptyfunc\colon \Moncat^{\op} \times \Moncat\to
\Moncat$ preserves weak equivalences between fibrant objects.

\end{thmlist}
\end{proposition}

\begin{proof}
Applying the Pushout Product Axiom to the map from the initial object
to a cofibrant object $X$ gives that $\monprod $ preserves
cofibrations and acyclic cofibrations in either variable when the
other variable is cofibrant.  A second application of the Pushout
Product Axiom then indicates that $\monprod $ preserves acyclic
cofibrations between cofibrant objects.  Coupled with this
observation, Brown's Lemma \cite[9.9]{modcat} implies that
$\monprod $ preserves weak equivalences between cofibrant objects.  A
similar argument, using the Enrichment Axiom \ENa, proves the claim
about $\emptyfunc$.
\end{proof} 

Parts (d) and (e) of this proposition indicate that $\monprod $ and
$\emptyfunc$ satisfy Quillen's criterion (Proposition~1 in
\cite[p. I.4.2]{quil} or Proposition~9.3 in \cite{modcat}) for the
existence of a left derived functor $\emptydermonprod$ and a right
derived functor $\emptyderfunc$, respectively.  Moreover, these
derived functors can be constructed so that
\[
\dermonprod{X}{Y} = X \monprod Y \tand \derfunc{Y}{Z} = \func{Y}{Z}
\]
whenever $X$ and $Y$ are cofibrant objects of $\Moncat$ and $Z$ is a
fibrant object of $\Moncat$.  In particular, if $Y$ is cofibrant, then
$\dermonprod{(-)}{Y}$ and $\derfunc{Y}{-}$ are the left and right
derived functors of $(-)\monprod Y$ and $\func{Y}{-}$ respectively.
Quillen's criterion for adjoint derived functors (Theorem~3 in
\cite[pp. I.4.5ff]{quil} or Theorem~9.7 in \cite{modcat}) then implies
that $\dermonprod{(-)}{Y}$ and $\derfunc{Y}{-}$ are adjoint.  These
observations are summarized in the following proposition.

\begin{proposition}\label{derivedbifunctors}
If $\Moncat$ is a monoidal model category, then the left derived
functor $\emptydermonprod$ of $\monprod$ and the right derived functor
$\emptyderfunc$ of $\emptyfunc$ exist and give a parametrized
adjunction
\[
\Ho\Moncat(\dermonprod{X}{Y},Z) \iso \Ho\Moncat(X,\derfunc{Y}{Z})
\]
\end{proposition}

This parametrized adjunction is the source of an easy proof of the
equivalence of the Unit Axioms \UF and \UC.  The functors $(-)\monprod
\cofunit$ and $\func{\cofunit}{-}$ are Quillen adjoints, and induce
the adjoint functors $\dermonprod{(-)}{\cofunit}$ and
$\derfunc{\cofunit}{-}$ on the homotopy category $\Ho\Moncat$.  The
map in axiom \UF is a natural transformation in $\Ho\Moncat$ from the
identity functor to $\derfunc{\cofunit}{-}$.  Axiom \UF asserts that
this natural transformation is a natural isomorphism.  The map in
axiom \UC is the adjoint natural transformation from
$\dermonprod{(-)}{\cofunit}$ to the identity.  That axiom asserts that
this adjoint natural transformation is a natural isomorphism.  Since
each natural transformation is a natural isomorphism if and only if
its adjoint is, the two axioms are equivalent.

From this it follows that the unit isomorphism for $\monprod$ induces
a unit isomorphism for $\emptydermonprod$.  Using the description of
the derived functor $\emptydermonprod$ in terms of cofibrant
approximations, it is straightforward to check that the associativity
isomorphism for $\monprod$ induces an associativity isomorphism for
$\emptydermonprod$.  Combined with
Proposition~\ref{derivedbifunctors}, these observations prove most of
the following result.  For a more complete discussion, see
\cite[4.3.2]{HvModCat}.

\begin{proposition}\label{HoMMndl}
Let $\Moncat$ be a monoidal model category.  The derived product
$\emptydermonprod$ and the derived function objects $\derfunc{-}{-}$
provide the homotopy category $\Ho\Moncat$ with a closed symmetric
monoidal structure.  Moreover, the localization functor
$\Moncat\to\Ho\Moncat$ is lax symmetric monoidal.
\end{proposition}


\section{Enriched model categories}\label{EnrModCatSec}

Although most of the main results stated in the introduction only make
sense for module categories, the most basic result,
Theorem~\ref{intenr}, applies more generally to closed model
categories enriched over a monoidal model category.  Moreover, since
the enrichment of the derived balanced product and function functors
described in our main results concisely encodes much of the coherence
among the derived functors we discuss, it is particularly convenient
to work in the context of enriched categories as much as possible.
Our first objective in this section is to introduce axioms for the
interaction of a model category structure with an enriched category
structure which imply the ``expected'' relationship between the
homotopy category, the enrichment, and the homotopy category of the
enriching monoidal model category.  We begin the discussion of this
relationship with Theorem~\ref{mainenr}, the generalization of
Theorem~\ref{intenr} to enriched model categories.  The discussion
then continues in the next section with Theorem~\ref{localenr}, which
states the universal property of the enrichment of the homotopy
category, and with a study of enrichments of derived functors.
Finally, we conclude the discussion in Section~\ref{BiFunSec} with a
study of enriched derived bifunctors and enriched
parametrized adjunctions.

Recall that a category $\EC$ enriched over a closed symmetric monoidal
category $\Moncat$ consists of:
\begin{hyplist}
\item A class $\ob{\EC}$ of \defterm{objects} of $\EC$,
\item For each $\objectone$, $\objecttwo$ in $\ob{\EC}$ a \defterm{mapping
object} $\efunc{\objectone}{\objecttwo}$ in $\Moncat$,
\item A \defterm{composition law} given by maps
\[
\circ \colon \efunc{\objecttwo}{\objectthree}\monprod
\efunc{\objectone}{\objecttwo} \to \efunc{\objectone}{\objectthree}
\]
in $\Moncat$ for each $\objectone$, $\objecttwo$, $\objectthree$ in
$\ob{\EC}$, and
\item \defterm{Identity morphisms}, which are maps $\id_{\objectone}\colon
\unit \to\efunc{\objectone}{\objectone}$ in $\Moncat$ for each
$\objectone$ in $\ob{\EC}$,
\end{hyplist}
These morphisms are required to satisfy the appropriate associativity
and unit conditions (see, for example, \cite[1.2]{BasicEnrch}).

The ordinary category underlying $\EC$ has the same objects as $\EC$
and morphism sets given by
\begin{equation*}
\EC(\objectone,\objecttwo)=\Moncat(\unit,\efunc{\objectone}{\objecttwo}).
\end{equation*}
The composition law and identity morphisms for this underlying
category are derived from the composition law and identity morphisms
in $\Moncat$ above.

More informally, an enrichment over $\Moncat$ of an ordinary category
$\EC$ is an isomorphism (or merely equivalence) between $\EC$ and the
underlying category of a $\Moncat$-enriched category.  For example,
$\Moncat$ is enriched over itself by the isomorphism
\[
\Moncat(\unit,\func{X}{Y})\iso \Moncat (\unit \monprod X,Y) \iso
\Moncat(X,Y).
\]

The following definition describes the standard procedure for
pushing enrichments forward along a monoidal functor.

\begin{definition}\label{EnrPFDef}
Let \map{\lambda}{\Moncat}{\Moncatb} be a lax symmetric monoidal
functor between two symmetric monoidal closed categories \Moncat and
$\Moncatb$.  Let $\EC$ be a category enriched over $\Moncat$.  The
induced category $\lambda_*\EC$ enriched over $\Moncatb$ has the same
object set as $\EC$ and morphism objects in $\Moncatb$ given by
\begin{equation*}
\efunc[{\lambda_*\EC}]{\objectone}{\objecttwo} =
\lambda\left(\efunc{\objectone}{\objecttwo}\right).
\end{equation*}
The composition and identity maps in \Moncatb for $\lambda_*\EC$ are
obtained by applying $\lambda$ to the analogous maps for $\EC$ in
\Moncat and composing with the appropriate morphisms giving $\lambda$
its monoidal structure.  There is a canonical functor from the
underlying category of $\EC$ to that of $\lambda_*\EC$ which is the
identity on objects and on morphisms is
\begin{eqnarray*}
\EC(\objectone,\objecttwo)=\Moncat(\unit_{\Moncat},\efunc{\objectone}{\objecttwo}) & \lrarrow[\lambda] & \Moncatb\left(\lambda\unit_{\Moncat},\lambda\left(\efunc{\objectone}{\objecttwo}\right)\right) \\ 
& \lrarrow &
\Moncatb\left(\unit_{\Moncatb},\lambda\left(\efunc{\objectone}{\objecttwo}\right)\right)=\lambda_*\EC(\objectone,\objecttwo),
\end{eqnarray*}
where the second map comes from the unit map for $\lambda$.
\end{definition}

The monoidal functor of interest to us is the localization functor
\map{\lambda}{\Moncat}{\Ho\Moncat} associated to a monoidal model
category \Moncat.  If $\EC$ is enriched over \Moncat, then
$\lambda_*\EC$ is a sort of homotopy category.  For example, when
\Moncat is the monoidal model category of spaces, the enrichment of
$\EC$ over \Moncat is given by function spaces
$\efunc{\objectone}{\objecttwo}$.  The morphism sets of $\lambda_*\EC$
are then the path components of these function spaces.  Thus, when
$\lambda$ is the localization functor for a monoidal model category,
$\lambda_*\EC$ is a natural generalization of the traditional
notion of a homotopy category.  If $\EC$ also carries a closed model
structure, then it is natural to inquire about the relationship
between $\Ho\EC$ and $\lambda_*\EC$.  Without some restrictions on the
model structure on $\EC$, there need not even be a functor comparing
$\Ho\EC$ and $\lambda_*\EC$.  However, there is an obvious
generalization of the Enrichment Axiom for monoidal model categories
to the context of closed model categories enriched over a monoidal
model category.  For the statement of this axiom, we need the
following generalization of the map $\boxfunc{g}{h}$ from Section
\ref{MndlModCatSec}.  Let \map{f}{A}{B}, \map{g}{K}{L}, and
\map{h}{X}{Y} be maps in a category $\EC$ enriched over a monoidal
model category \Moncat.  Then
\begin{equation*}
\boxefunc{g}{h}\colon \efunc{L}{X}\to
\efunc{K}{X}\times_{\efunc{K}{Y}}\efunc{L}{Y}
\end{equation*}
is the map defined by the pullback analogous to \eqref{diaboxfunc}
with \emptyefunc in place of \emptyfunc.  The generalization of the
Enrichment Axiom for monoidal model categories to the context of
enriched categories is

\begin{description}
\item[\ENCa] If \map{g}{K}{L} is a cofibration in $\EC$ and
\map{h}{X}{Y} is a fibration in $\EC$, then $\boxefunc{g}{h}$ is a
fibration in $\Moncat$.  Moreover, if either $g$ or $h$ is also a weak
equivalence, then so is $\boxefunc{g}{h}$.
\end{description}

The following analog of Proposition \ref{pppweq} describes the
implications of this Enrichment Axiom for the functor $\efunc{-}{-}$.
(Recall for part~(c) that a fibrant object in
$\Moncat^{\op}$ is a cofibrant object in $\Moncat$.)

\begin{proposition}\label{EnrImp}
Let $\EC$ be a closed model category enriched over a monoidal model
category $\Moncat$ such that the Enrichment Axiom is satisfied.
\begin{thmlist}
\item If $\objectone$ is cofibrant, then the functor
$\efunc{\objectone}{-}$ preserves fibrations and acyclic fibrations.

\item If $\objecttwo$ is fibrant, then the functor
$\efunc{-}{\objecttwo}$ converts cofibrations and acyclic cofibrations
in $\EC$ into fibrations and acyclic fibrations in $\Moncat$,
respectively.

\item\label{eFunweProp} The functor $\efunc{-}{-}\colon \EC^{\op}
\times \EC\to \Moncat$ preserves weak equivalences between fibrant
objects.

\end{thmlist}
\end{proposition}

Proposition \ref{EnrImp}(\ref{eFunweProp}) indicates that the restriction of the
composite functor
\begin{equation*}
\EC^{\op} \times \EC \lrarrow[\efunc{-}{-}] \Moncat \lrarrow[\lambda]
\Ho \Moncat
\end{equation*}
to the full subcategory of $\EC_{cf}^{\op} \times \EC_{cf}$ consisting
of pairs $(C,D)$ such that $C$ and $D$ are both cofibrant and fibrant
in $\EC$ converts weak equivalences into isomorphisms.  It follows
from the universal property of localization that the functor $\EC_{cf}
\lrarrow \lambda_*\EC_{cf}$ factors through the category $\Ho\EC_{cf}$
to give a comparison functor
\begin{equation*}
\map{\Upsilon}{\Ho\EC_{cf}}{\lambda_*\EC_{cf}}.
\end{equation*}
This comparison functor is the subject of the following \defterm{Homotopy/Unit
Axiom}:
\begin{description}
\item[\UE] The functor \map{\Upsilon}{\Ho\EC_{cf}}{\lambda_*\EC_{cf}}
is an isomorphism of categories.
\end{description}
In other words, the Homotopy/Unit Axiom requires that whenever $C$ and
$D$ are cofibrant-fibrant objects of $\EC$, the map $\Ho\EC(C,D)\to
\lambda_{*}\EC(C,D) $ is a bijection.

This axiom turns out to generalize the Unit Axiom in the definition of
a monoidal model category.  It is shown below that it is equivalent to
both of the more obvious generalizations of the Unit Axiom that become
available when $\EC$ is tensored or cotensored over \Moncat.
Together, the Enrichment Axiom and Homotopy/Unit Axiom suffice to
describe the model structures on enriched categories which give
homotopy categories that appropriately preserve the enrichment.

\begin{definition}\label{defenr}
Let $\Moncat$ be a monoidal model category and let $\EC$ be a closed
model category that is enriched over $\Moncat$.  Then $\EC$ is an
\defterm{enriched model category} if it satisfies the Enrichment Axiom \ENCa and
the Homotopy/Unit Axiom \UE.
\end{definition}

When $\EC$ has tensors or cotensors, both the Enrichment Axiom \ENCa
and the Homotopy/Unit Axiom \UE have alternate forms that are easier
to verify in practice.  For the statements of these alternative forms,
we need the following generalizations of the maps $\boxfunc{g}{h}$ and
$f\boxmonprod g$ defined in Section \ref{MndlModCatSec}. Let
\map{f}{A}{B} and \map{h}{X}{Y} be maps in $\EC$ and \map{g}{K}{L} be
a map in $\Moncat$.  Then the maps
\begin{equation*}
\boxcot{g}{h}\colon \coten{L}{X}\to
\coten{K}{X}\times_{\coten{K}{Y}}\coten{L}{Y}
\end{equation*}
and
\begin{equation*}
\boxten{g}{f} \colon (\ten{L}{A}) \cup_{\ten{K}{A}}(\ten{K}{B})\to
\ten{L}{B}
\end{equation*}
are defined as the pullback analogous to \eqref{diaboxfunc} and as the
pushout analogous to \eqref{diaboxprod} (with $\emptytensor$ replacing
$\monprod$), respectively. 

The following proposition provides the alternative forms of the
Enrichment Axiom.  It follows easily from the characterization of
(acyclic) cofibrations and (acyclic) fibrations in $\Moncat$ in terms
of lifting properties, using the tensor or cotensor adjunction.

\begin{proposition}\label{altenrEN}
Let $\Moncat$ be a monoidal model category and let $\EC$ be a closed
model category that is also enriched over $\Moncat$.
\begin{thmlist}

\item If $\EC$ has tensors, then the Enrichment Axiom \ENCa is
equivalent to the following Pushout Tensor Product Axiom:

\begin{description}
\item[\tPPa]If $f\colon A\to B$ is a cofibration in $\EC$ and $g\colon
K\to L$ is a cofibration in $\Moncat$, then \boxten{g}{f} is a
cofibration in $\EC$. Moreover, if either $f$ or $g$ is also a weak
equivalence, then so is \boxten{g}{f}.
\end{description}

\item If $\EC$ has cotensors, then the Enrichment Axiom \ENCa is
equivalent to the following Cotensor Axiom:

\begin{description}
\item[\cSMa] If \map{g}{K}{L} is a cofibration in $\Moncat$ and
\map{h}{X}{Y} is a fibration in $\EC$, then \boxcot{g}{h} is a
fibration in $\EC$. Moreover, if either $g$ or $h$ is also a weak
equivalence, then so is \boxcot{g}{h}.
\end{description}

\end{thmlist}
\end{proposition}

The following proposition provides the alternative forms of the
Homotopy/Unit Axiom \UE:

\begin{proposition}\label{altenrUE}
Let $\Moncat$ be a monoidal model category and let $\EC$ be a closed
model category that is enriched over $\Moncat$ and satisfies the
Enrichment Axiom \ENCa.
\begin{thmlist}

\item If $\EC$ has tensors, then the Homotopy/Unit Axiom is equivalent
to the following Tensor Unit Axiom
\begin{description}
\item[\tUC]When $C$ is cofibrant in $\EC$, the map $\cum\colon
\tensor{C}{\cofunit}\to \tensor{C}{\unit}\iso C$ is a weak
equivalence.
\end{description}

\item If $\EC$ has cotensors, then the Homotopy/Unit Axiom is
equivalent to the following Cotensor Unit Axiom

\begin{description}
\item[\cUF]When $D$ is fibrant in $\EC$, the map $\cuf\colon D \iso
\cotensor{\unit}{D} \to\cotensor{\cofunit}{D}$ is a weak equivalence.
\end{description}

\end{thmlist}
\end{proposition}

The proof of this result is similar to the proof of the equivalence of
the two unit axioms \UF and \UC given in Section \ref{MndlModCatSec}.
It makes use of the adjunction relating the derived functors of the
tensor and cotensor functors, and this proof is delayed until after
our discussion of the existence of these derived functors.  The
following extension of Proposition~\ref{pppweq} to the context of
tensors and cotensors is needed in the discussion.  Its proof again
follows by adjunction from the characterization of (acyclic)
cofibrations and (acyclic) fibrations in terms of lifting.
(Recall for part~(b) that a fibrant object in $\Moncat^\op$ is a
cofibrant object of $\Moncat$.)

\begin{proposition}\label{TenCotImp}
Let $\Moncat$ be a monoidal model category and let $\EC$ be a closed
model category that is enriched over $\Moncat$ and satisfies the
Enrichment Axiom \ENCa.
\begin{thmlist}
\item\label{TenQAdj} Assume $\EC$ has tensors.  If $\objectone$ and
$X$ are cofibrant in $\EC$ and $\Moncat$, respectively, then
$\tensor{\objectone}{(-)}$ and $\tensor{(-)}{X}$ preserve cofibrations
and acyclic cofibrations.  Also, the functor $\emptytensor\colon
\EC\times \Moncat\to\EC$ preserves weak equivalences between cofibrant
objects.

\item Assume $\EC$ has cotensors. If $X$ is cofibrant in $\Moncat$,
then $\cotensor{X}{-}$ preserves fibrations and acyclic fibrations in
$\EC$.  If $\objecttwo$ is fibrant in $\EC$, then
$\cotensor{-}{\objecttwo}$ converts cofibrations and acyclic
cofibrations in $\Moncat$ into fibrations and acyclic fibrations in
$\EC$, respectively.  Also, the functor $\emptycotensor\colon
\Moncat^{\op}\times \EC\to \EC$ preserves weak equivalences between
fibrant objects.

\end{thmlist}
\end{proposition}

Proposition \ref{TenCotImp} implies that the tensor and cotensor
adjunctions for $\EC$ are Quillen adjunctions in each variable.  This
observation and an argument analogous to the proof of
Proposition~\ref{derivedbifunctors} proves the following proposition.

\begin{proposition}\label{tensoradj}
Let $\Moncat$ be a monoidal model category and let $\EC$ be a closed
model category that is enriched over $\Moncat$ and satisfies the
Enrichment Axiom \ENCa.
\begin{thmlist}
\item If $\EC$ has tensors then the left derived functor
$\emptydertensor$ of $\emptytensor$ exists and with $\emptyderefunc$
gives a parametrized adjunction
\[
\Ho\EC(\dertensor{\objectone}{X},\objecttwo)\iso
\Ho\Moncat(X,\derefunc{\objectone}{\objecttwo}).
\]
\item If $\EC$ has cotensors then the right derived functor
$\emptydercotensor$ of $\emptycotensor$ exists and with
$\emptyderefunc$ gives a parametrized adjunction
\[
\Ho\EC(\objectone,\dercotensor{X}{\objecttwo})\iso
\Ho\Moncat(X,\derefunc{\objectone}{\objecttwo}).
\]
\end{thmlist}
\end{proposition}

A cofibrant approximation \vmap{\cofunit}{\unit} to the unit $\unit$
for \Moncat yields natural transformations
\begin{equation*}
\cum\colon \tensor{C}{\cofunit}\to \tensor{C}{\unit}\iso C
\end{equation*}
and
\begin{equation*}
\cuf\colon D \iso \cotensor{\unit}{D} \to\cotensor{\cofunit}{D}
\end{equation*}
relating the inclusion functor on $\Ho\EC_{cf}\to \Ho\EC$ to the functors
$\tensor{-}{\cofunit}$ and $\cotensor{\cofunit}{-}$.  The axioms \UC
and \UF assert that these natural transformations are natural
isomorphisms.  The adjunction of Proposition \ref{tensoradj} allows us
to relate these natural transformations to the comparison functor
\map{\Upsilon}{\Ho\EC_{cf}}{\lambda_*\EC_{cf}} via the following
commuting diagram:
\begin{equation}\label{UpsDiag}
\begin{gathered}
\xymatrix@=1em{ {\Ho\EC_{cf}(C,D) \ar[rr]^-{(\cuf)_*} \ar[dd]_{\cum^*}
\ar[dr]^-{\Upsilon}} & {}
	& {\Ho\EC_{cf}(C,\cotensor{\cofunit}{D}) \ar[dd]^{\iso}} \\
{} & {\lambda_*\EC(C,D) \ar[dr]^-{=} }
	& {} \\	
{\Ho\EC(\tensor{C}{\cofunit},D) \ar[rr]_-{\iso}} & {} &
{\Ho\Moncat(\cofunit,\derefunc{C}{D})} }
\end{gathered}
\end{equation}
(Only the relevant part of this diagram exists when $\EC$ has tensors
but not cotensors or vice-versa.)  Clearly each of the maps
$(\cuf)_*$, $\Upsilon$, and $\cum^*$ in this diagram is an isomorphism
if and only if the either of the other
maps is also an isomorphism.  This implies that \UE is equivalent to
\cUF and \tUC whenever either axiom makes sense.

Our motivating examples of enriched model categories are provided by
the following result.

\begin{proposition}\label{leftmodenr}
Let $\Moncat$ be a monoidal model category, let $\monoid$ be a monoid
in $\Moncat$, and assume the category $\Leftmod$ of left $\monoid$
modules is a closed model category with fibrations and weak
equivalences created in $\Moncat$.  Then $\Leftmod$ is an enriched
model category.
\end{proposition}

\begin{proof}
$\Leftmod$ has tensors induced by $\monprod$ and cotensors induced by
$\emptyfunc$.  Moreover, the Enrichment Axiom \ENa for $\Moncat$
implies the Cotensor Axiom \cSMa for $\Leftmod$ and the Unit Axiom \UF
for $\Moncat$ implies the Cotensor Unit Axiom \cUF for $\Leftmod$.
\end{proof}

The following is our fundamental result about enriched model
categories.  Coupled with the previous proposition, it implies
Theorem~\ref{intenr} of the introduction.

\begin{theorem}\label{mainenr}
Let $\EC$ be an enriched model category over a monoidal model category
$\Moncat$.  Then the right derived functor $\emptyderefunc$ of
$\emptyefunc$ exists and enriches $\Ho\EC$ over $\Ho\Moncat$.
Further, if $\EC$ is tensored or cotensored over $\Moncat$, then
$\Ho\EC$ is likewise tensored or cotensored over $\Ho\Moncat$.
\end{theorem}

\begin{proof}
It follows from Proposition \ref{eFunweProp} that the right derived
functor $\emptyderefunc$ of $\emptyefunc$ exists and $\derefunc{X}{Y}$
may be computed as $\efunc{QX}{RY}$ where $QX\to X$ is a cofibrant
approximation and $Y\to RY$ is a fibrant approximation.  Every object
$\objectone$ of $\Ho\EC$ is isomorphic in $\Ho\EC$ to an object
$RQ\objectone$ of $\Ho\EC_{cf}$; conjugating by these isomorphisms
gives both an isomorphism $\derefunc{RQ\objectone}{RQ\objecttwo}\iso
\derefunc{\objectone}{\objecttwo}$ and Quillen's equivalence of
$\Ho\EC_{cf}$ with $\Ho\EC$.  The Homotopy/Unit Axiom requires that
\map{\Upsilon}{\Ho\EC_{cf}}{\lambda_*\EC_{cf}} is an isomorphism, and
viewing $\lambda_{*}\EC_{cf}$ as a full enriched subcategory of
$\lambda_{*}\EC$, it follows that $\Ho\EC_{cf}$ and therefore $\Ho\EC$
is enriched by $\emptyderefunc$.  This completes the proof of the
first statement of the theorem.  Proposition \ref{tensoradj} is a
first step toward proving that $\Ho\EC$ is tensored and/or cotensored
over $\Ho\Moncat$ when $\EC$ is.  However, the adjunctions provided by
that proposition are ordinary, rather than enriched adjunctions.  To
complete the proof of the theorem, we prove a stronger version of
Proposition \ref{tensoradj} with enriched adjunctions as
Corollary~\ref{mainenrtwo} in our discussion of enriched functors and
adjunctions in the next section.
\end{proof}

For concreteness and to introduce notation used in the next section,
we describe in more detail the composition law constructed in the
previous proof.  For each object $\objectone$ of $\EC$, choose and fix
an acyclic fibration $q_{\objectone}\colon Q\objectone \to \objectone$
with $Q\objectone$ cofibrant (with $q_{\objectone}$ the identity if
$\objectone$ is cofibrant), an acyclic cofibration
$r_{\objectone}\colon \objectone\to R\objectone$ with $R\objectone$
fibrant (with $r_{\objectone}$ the identity if $\objectone$ is
fibrant), and a factorization $s_{\objectone}\colon RQ\objectone \to
R\objectone$ of the composite $Q\objectone\to R\objectone$, i.e., a
map $s_{\objectone}$ making the diagram on the left commute.
\[ \xymatrix{%
Q\objectone\ar@{->>}[d]_{q_{\objectone}}^{\sim}
\ar@{>->}[r]^{r_{Q\objectone}}_{\sim}
&RQ\objectone\ar@{-->}[d]^{s_{\objectone}}_{\sim}&&
Q\objectone\ar@{>->}[d]_{r_{Q\objectone}}^{\sim}
\ar[r]^{r_{\objectone}\circ q_{\objectone}}&R\objectone\ar@{->>}[d]\\
\objectone\ar@{>->}[r]_{r_{\objectone}}^{\sim}&R\objectone&&
RQ\objectone\ar[r]\ar@{-->}[ur]^{s_{\objectone}}&\relax* } 
\]
Such a factorization exists by the lifting property of cofibrations with
respect to acyclic fibrations illustrated on the diagram on the right
above. Note that $s_{\objectone}$ is a weak equivalence by the
two-out-of-three property. We choose $s_{\objectone}$ to be the
identity when $\objectone$ is cofibrant.  Then
$s_{\objectone}^{-1}\circ r_{\objectone} =r_{Q\objectone}\circ
q_{\objectone}^{-1}$ is an isomorphism in $\Ho\EC$ from $\objectone$
to the cofibrant-fibrant object $RQ\objectone$.

The purpose of the choice of the maps $s$ is that it allows us to
identify the isomorphism
\[
(r_{Q\objectone}\circ q_{\objectone}^{-1})^{*} (r_{Q\objecttwo}\circ
q_{\objecttwo}^{-1})^{-1}_{*} \colon
\derefunc{RQ\objectone}{RQ\objecttwo}\to
\derefunc{\objectone}{\objecttwo}
\]
in the proof of Theorem~\ref{mainenr} above with the map
\[
\derefunc{RQ\objectone}{RQ\objecttwo}=\efunc{RQ\objectone}{RQ\objecttwo}
\lrarrow[(s_{\objecttwo})_{*}r_{Q\objectone}^{*}]
\efunc{Q\objectone}{R\objecttwo}=\derefunc{\objectone}{\objecttwo}.
\]
The composition in $\Ho\EC$ therefore fits into the following
commutative diagram in $\Ho\Moncat$, where the dotted arrows are the
inverses of the isomorphisms indicated by the corresponding backward
solid arrows.
\begin{equation*}
\begin{gathered}
\def\objectstyle{\scriptstyle}\def\labelstyle{\scriptscriptstyle}
\xymatrix{%
\relax\dermonprod{\derefunc{\objecttwo}{\objectthree}}{\derefunc{\objectone}{\objecttwo}}
\ar@{=}[d]\ar@{-->}[rr]^{\circ}
&&\derefunc{\objectone}{\objectthree}\ar@{=}[d]
\\
\relax\dermonprod{\efunc{Q\objecttwo}{R\objectthree}}{\efunc{Q\objectone}{R\objecttwo}}
\backisoabove[r]{\dermonprod{r_{Q\objecttwo}^{*}}{(s_{\objecttwo})_{*}}}
\backisobelow[dr]{\dermonprod{r_{Q\objecttwo}^{*}(s_{\objectthree})_{*}}{(s_{\objecttwo})_{*}r_{Q\objectone}^{*}}\qquad}
&\relax\dermonprod{\efunc{RQ\objecttwo}{R\objectthree}}{\efunc{Q\objectone}{RQ\objecttwo}}
\ar[r]^(.60){\circ}
\backisoabove[d]{\dermonprod{(s_{\objectthree})_{*}}{r_{Q\objectone}^{*}}}
&\relax\efunc{Q\objectone}{R\objectthree}
\backisoabove[d]{{(s_{\objectthree})_{*}}{r_{Q\objectone}^{*}}}
\\
&
\relax\dermonprod{\efunc{RQ\objecttwo}{RQ\objectthree}}{\efunc{RQ\objectone}{RQ\objecttwo}}
\ar[r]_(.60){\circ} &\relax\efunc{RQ\objectone}{RQ\objectthree}
} \end{gathered}
\end{equation*}
The horizontal arrows in the bottom right square are the composition
in $\lambda_{*}\EC$.  We can regard the middle row as a definition of
the composition in the enrichment of $\Ho\EC$.


\section{Enriched functors and enriched derived
functors}\label{EnrFunSec}

This section continues the study of enriched model categories with a
discussion of enriched functors.  We characterize the enrichment of
the homotopy category of an enriched model category in terms of a
universal property with respect to enriched functors.  This leads to a
generalization to enriched functors of Quillen's criterion for the
existence of derived functors and a corresponding theory of enriched
Quillen adjunctions.

Recall that, for categories $\EC$ and $\ED$ enriched over $\Moncat$,
an \defterm{enriched functor} $\Phi \colon \EC\to\ED$ consists of a function
$\Phi \colon \ob{\EC}\to \ob{\ED}$ together with maps
\[
\Phi_{\objectone,\objectone'}\colon \efunc{\objectone}{\objectone'}\to
\efunc[\ED]{\Phi \objectone}{\Phi \objectone'}
\]
in $\Moncat$ consistent with the identity morphisms and composition
law.  We also write $\Phi$ for the functor on the underlying
categories; this underlying functor is given by
$\Moncat(\unit,\Phi_{\objectone,\objectone'})$.  More generally, when
$\EC$ is enriched over $\Moncat$, $\ED$ is enriched over $\Moncatb$,
and $\lambda \colon \Moncat\to\Moncatb$ is a lax symmetric monoidal
functor, a $\lambda$-enriched functor $\Phi \colon \EC\to\ED$ (or
\defterm{$\Moncatb$-enriched}, when $\lambda$ is understood) consists of a
function $\Phi$ on objects and maps in $\Moncatb$
\[
\Phi_{\objectone,\objectone'}\colon \lambda
(\efunc{\objectone}{\objectone'})\to \efunc[\ED]{\Phi \objectone}{\Phi
\objectone'}
\]
consistent with the identity morphisms and composition law.  The
following well-known proposition essentially provides an equivalent
alternate definition of a $\lambda$-enriched functor in terms of the
$\Moncatb$-enriched category $\lambda_{*}\EC$ of the previous section.

\begin{proposition}
For any lax symmetric monoidal functor $\lambda \colon
\Moncat\to\Moncatb$ and any $\Moncat$-enriched category $\EC$,
$\lambda$ induces a $\lambda$-enriched functor $\EC\to
\lambda_{*}\EC$, and this $\lambda$-enriched functor is initial.  In
other words, for any $\Moncatb$-enriched category $\ED$, any
$\lambda$-enriched functor $\EC\to \ED$ factors uniquely through an
$\Moncatb$-enriched functor $\lambda_{*}\EC\to \ED$.
\end{proposition}

We are mainly concerned with the case where $\lambda$ is the
localization functor $\Moncat\to\Ho\Moncat$.  Using this special case
of a $\lambda$-enriched functor, we can identify the homotopy category
of an enriched model category by a universal property.  To avoid
confusion with the localization functor $\lambda \colon \Moncat\to
\Ho\Moncat$, we denote the localization functor $\EC\to \Ho\EC$ as
$\gamma$.

\begin{theorem}\label{localenr}
Let $\EC$ be an enriched model category over a monoidal model category
$\Moncat$.  The localization functor $\gamma \colon \EC\to\Ho\EC$ is
$\lambda$-enriched and is the initial $\lambda$-enriched functor that
sends weak equivalences to isomorphisms. In other words, for any
$\Ho\Moncat$-enriched category $\EH$, any $\lambda$-enriched functor
$\EC\to\EH$ that sends weak equivalences to isomorphisms factors
uniquely through a $\Ho\Moncat$-enriched functor $\Ho\EC\to \EH$.
\end{theorem}

\begin{proof}
The enriched localization functor is given by the universal maps
$\efunc{\objectone}{\objecttwo}\to \derefunc{\objectone}{\objecttwo}$
of the right derived functor; we need to check that these maps
assemble into an enriched functor.  The fact that they preserve the
identity morphisms is clear, and so it suffices to check that they preserve
composition.  Consider the following diagram in
$\Ho\Moncat$ written in the notation introduced at the end of the
previous section.
\begin{equation}\label{enrcompdiag}
\def\objectstyle{\scriptstyle}\def\labelstyle{\scriptscriptstyle}
\begin{gathered}
\xymatrix@C+1em{%
\relax\dermonprod{\efunc{\objecttwo}{\objectthree}}{\efunc{\objectone}{\objecttwo}}
\ar[rr]^{\circ}\ar[d]_{\dermonprod{(r_{\objectthree})_{*}}{q_{\objectone}^{*}}}
&&
\relax\efunc{\objectone}{\objectthree}\ar[dd]^{{(r_{\objectthree})_{*}}{q_{\objectone}^{*}}}
\\
\relax\dermonprod{\efunc{\objecttwo}{R\objectthree}}{\efunc{Q\objectone}{\objecttwo}}
\ar@/^4em/[drr]^{\circ}
\ar[d]_{\dermonprod{q_{\objecttwo}^{*}}{(r_{\objecttwo})_{*}}}
\backisobelow[r]{\dermonprod{r_{\objecttwo}^{*}}{(q_{\objecttwo})_{*}}}
&\relax\dermonprod{\efunc{R\objecttwo}{R\objectthree}}{\efunc{Q\objectone}{Q\objecttwo}}
\ar[d]^{\dermonprod{s_{\objecttwo}^{*}}{(r_{Q\objecttwo})_{*}}}
\\
\relax\dermonprod{\efunc{Q\objecttwo}{R\objectthree}}{\efunc{Q\objectone}{R\objecttwo}}
\backisoabove[r]{\dermonprod{r_{Q\objecttwo}^{*}}{(s_{\objecttwo})_{*}}}
&\relax\dermonprod{\efunc{RQ\objecttwo}{R\objectthree}}{\efunc{Q\objectone}{RQ\objecttwo}}
\ar[r]^(.60){\circ} &\relax\efunc{Q\objectone}{R\objectthree} }
\end{gathered}
\end{equation} 
The top row is the composition in $\lambda_{*}\EC$ and the bottom row
is essentially the composition in $\Ho\EC$.  The square
\[
\xymatrix@C+1em{%
\relax\dermonprod{\efunc{R\objecttwo}{R\objectthree}}{\efunc{Q\objectone}{Q\objecttwo}}
\ar[d]_{\dermonprod{s_{\objecttwo}^{*}}{(r_{Q\objecttwo})_{*}}}
&\relax\dermonprod{\efunc{\objecttwo}{R\objectthree}}{\efunc{Q\objectone}{\objecttwo}}
\ar[d]^{\circ}
\backisobelow[l]{\dermonprod{r_{\objecttwo}^{*}}{(q_{\objecttwo})_{*}}}
\\
\relax\dermonprod{\efunc{RQ\objecttwo}{R\objectthree}}{\efunc{Q\objectone}{RQ\objecttwo}}
\ar[r]_(.60){\circ} &\relax\efunc{Q\objectone}{R\objectthree} } \]
(where the right vertical arrow is the curved arrow in
diagram~\eqref{enrcompdiag}) commutes by dinaturality since
$r_{\objecttwo}\circ q_{\objecttwo}=s_{\objecttwo}\circ
r_{Q\objecttwo}$, and the remaining squares commute by naturality.  It
follows that $\EC\to \Ho\EC$ is $\Ho\Moncat$-enriched.

Given any $\lambda$-enriched functor $\Phi \colon \EC\to \EH$ that
sends weak equivalences to isomorphisms, the natural maps
$\Phi(\dermonprod{(r_{\objectthree})_{*}}{q_{\objectone}^{*}})$, $\Phi
(\dermonprod{q_{\objecttwo}^{*}}{(r_{\objecttwo})_{*}})$, and $\Phi
((r_{\objectthree})_{*}q_{\objectone}^{*})$ are isomorphisms in
$\Ho\Moncat$, and it follows from diagram~\eqref{enrcompdiag} that
$\Phi$ factors uniquely through a $\Ho\Moncat$-enriched functor
$\Ho\EC\to\EH$.
\end{proof}

Next we discuss derived functors in the enriched model category
context.  We concentrate our discussion on left derived functors to
avoid tedious repetition.  Recall that for a functor $\Phi \colon
\EC\to \EH$, the left derived functor $L\Phi\colon \Ho\EC\to\EH$ (if
it exists) is defined is to be the right Kan extension of $\Phi$ along
the localization functor $\gamma \colon \EC\to \Ho\EC$.  In other
words, the left derived functor (if it exists) as part of its
structure comes with a natural transformation $\phi \colon L\Phi\circ
\gamma \to \Phi$ which is final among natural transformations $F\circ
\gamma\to \Phi$.  The definition of enriched derived functors
therefore first requires review of the definition of enriched natural
transformations.

\begin{definition}
An \defterm{enriched natural transformation} $\alpha$ between enriched functors
$\Phi,\Phi'\colon \EC\to\ED$ is a natural transformation between the
underlying functors that makes the following diagram commute:
\[
\xymatrix{ {\efunc[{\EC}]{C}{C'} \ar[r]^-{\Phi_{C,C'}}
\ar[d]_-{\Phi'_{C,C'}} }
& {\efunc[{\ED}]{\Phi C}{\Phi C'} \ar[d]^-{(\alpha_{C'})_*}} \\
{\efunc[{\ED}]{\Phi' C}{\Phi' C'} \ar[r]_-{(\alpha_{C})^*}} &
{\efunc[{\ED}]{\Phi C}{\Phi' C'}} }
\]
If, instead, $\Phi$ and $\Phi'$ are $\lambda$-enriched functors, then
a \defterm{$\lambda$-enriched natural transformation} is a $\Moncatb$-enriched
natural transformation from $\Phi$ to $\Phi'$, considered as
$\Moncatb$-enriched functors out of $\lambda_*\EC$
\end{definition}

We offer the following definition in analogy with the definition of
left derived functor.

\begin{definition}
Let $\EC$ be an enriched model category over a monoidal model category
$\Moncat$.  Let $\EH$ be a $\Ho\Moncat$-enriched category and let
$\Phi \colon \EC\to \EH$ be a $\lambda$-enriched functor.  We say that
a $\Ho\Moncat$-enriched functor $L^{\Moncat}\Phi\colon \Ho\EC\to \EH$
and $\lambda$-enriched natural transformation $\phi^{\Moncat}\colon
L^{\Moncat}\circ \gamma \to \Phi$ forms the \defterm{enriched left derived
functor} of $\Phi$ when $\phi^{\Moncat}$ is final among
$\lambda$-enriched natural transformations $F\circ \gamma \to \Phi$.
In other words, given any $\Ho\Moncat$-enriched functor $F\colon
\Ho\EC\to \EH$ and $\lambda$-natural transformation $\alpha \colon
F\circ \gamma \to \Phi$, there exists a unique $\Ho\Moncat$-enriched
natural transformation $\theta \colon F\to L^{\Moncat}\Phi$ such that
$\phi^{\Moncat}\circ \theta =\alpha$.
\end{definition}

The enriched right derived functor is defined analogously, or
equivalently, as $R^{\Moncat}\Phi=(L^{\Moncat}\Phi^{\op})^{\op}$, for
$\Phi^{\op}\colon \EC^{\op}\to \EH^{\op}$.  Note that without further
hypotheses on $\Phi$, the underlying functor and natural
transformation of the enriched left derived functor need not agree
with the left derived functor of $\Phi$ when both exist.  In the case
when they do agree, we say that $L^{\Moncat}\Phi, \phi^{\Moncat}$
provide an \defterm{enrichment} of the derived functor.  Next we extend
Quillen's criterion for the existence of left derived functors to the
enriched context, and show that under its hypotheses, the enriched
left derived functor exists and provides an enrichment for the derived
functor.

Quillen's criterion for the existence of a left derived functor
asserts that when $\Phi \colon \EC\to\EH$ preserves weak equivalences
between cofibrant objects, the left derived functor exists and can be
computed using the cofibrant approximations $Q\objectone$.  In detail,
for each map $f\colon \objectone\to \objecttwo$ in $\EC$, we choose
$Qf\colon Q\objectone\to Q\objecttwo$ to be a lift of $f\circ
q_{\objectone}$, i.e., choose a function $Q_{\objectone,\objecttwo}$
making the following diagram commute.
\[
\xymatrix{%
\EC(\objectone,\objecttwo)
\ar@{-->}[rr]^{Q_{\objectone,\objecttwo}}
\ar[dr]_{q_{\objectone}^{*}}&&
\EC(Q\objectone,Q\objecttwo)\ar[dl]^{(q_{\objecttwo})_{*}}\\
&\EC(Q\objectone,\objecttwo) }
\]
Although $Q$ is not a functor, implicit in the statement and explicit
in the proof of Quillen's criterion is that when $\Phi$ preserves weak
equivalences between cofibrant objects, the composite $\Phi \circ Q$
becomes a functor and $\Phi(q)$ a natural transformation.  In the
enriched context, the map $(q_{\objecttwo})_{*}\colon
\efunc{Q\objectone}{Q\objecttwo}\to \efunc{Q\objectone}{\objecttwo}$ is
an acyclic fibration, and so is an isomorphism in $\Ho\Moncat$.  Thus,
there exists a unique map $Q_{\objectone,\objecttwo}$ in $\Ho\Moncat$
making the following diagram in $\Ho\Moncat$ commute.
\[
\xymatrix{%
\efunc{\objectone}{\objecttwo}
\ar@{-->}[rr]^{Q_{\objectone,\objecttwo}}
\ar[dr]_{q_{\objectone}^{*}}&&
\efunc{Q\objectone}{Q\objecttwo}\ar[dl]^{(q_{\objecttwo})_{*}}_{\simeq}\\
&\efunc{Q\objectone}{\objecttwo} }
\]
This leads to the following observation.

\begin{lemma}\label{Qenr}
There is an enriched functor $Q\colon \lambda_{*}\EC\to
\lambda_{*}\EC_{c}$ extending the function $Q$ on objects.  The maps
$q$ assemble to an enriched natural transformation from $Q$ to the
identity in $\lambda_{*}\EC$.
\end{lemma}

\begin{proof}
As indicated above, the enriched functor $Q$ is defined as the map in
$\Ho\Moncat$
\[ \xymatrix{%
Q_{\objectone,\objecttwo}\colon
\efunc{\objectone}{\objecttwo}\ar[r]^{q_{\objectone}^{*}}
&\efunc{Q\objectone}{\objecttwo}
\backisoabove[r]{(q_{\objecttwo})_{*}}
&\efunc{Q\objectone}{Q\objecttwo} } \] and it is clear from this
definition that $q$ is an enriched natural transformation provided
that $Q$ is an enriched functor.  To see that $Q$ is an enriched
functor, consider the following diagram in $\Ho\Moncat$.
\[ \def\objectstyle{\scriptstyle}\def\labelstyle{\scriptscriptstyle}
\xymatrix@C+1em{%
\relax\dermonprod{\efunc{\objecttwo}{\objectthree}}{\efunc{\objectone}{\objecttwo}}
\ar[rr]^{\circ}\ar[d]_{\dermonprod{\id}{q_{\objectone}^{*}}} &&
\relax\efunc{\objectone}{\objectthree}\ar[dd]^{q_{\objectone}^{*}}
\\
\relax\dermonprod{\efunc{\objecttwo}{\objectthree}}{\efunc{Q\objectone}{\objecttwo}}
\ar@/^4em/[drr]^{\circ} \ar[d]_{\dermonprod{q_{\objecttwo}^{*}}{\id}}
\backisobelow[r]{\dermonprod{\id}{(q_{\objecttwo})_{*}}}
&\relax\dermonprod{\efunc{\objecttwo}{\objectthree}}{\efunc{Q\objectone}{Q\objecttwo}}
\ar[d]^{\dermonprod{q_{\objecttwo}^{*}}{\id}}
\\
\relax\dermonprod{\efunc{Q\objecttwo}{\objectthree}}{\efunc{Q\objectone}{\objecttwo}}
\backisoabove[r]{\dermonprod{\id}{(q_{\objecttwo})_{*}}}
\backisobelow[d]{\dermonprod{(q_{\objecttwo})_{*}}{(q_{\objectthree})_{*}}}
&\relax\dermonprod{\efunc{Q\objecttwo}{\objectthree}}{\efunc{Q\objectone}{Q\objecttwo}}
\backisoabove[dl]{\dermonprod{(q_{\objectthree})_{*}}{\id}}
\ar[r]^(.60){\circ} &\relax\efunc{Q\objectone}{\objectthree}
\backisoabove[d]{\dermonprod{(q_{\objectthree})_{*}}{\id}}
\\
\relax\dermonprod{\efunc{Q\objecttwo}{Q\objectthree}}{\efunc{Q\objectone}{Q\objecttwo}}
\ar[rr]_{\circ} && \relax\efunc{Q\objectone}{Q\objectthree} }
\]
This diagram, like diagram~\eqref{enrcompdiag}, commutes by naturality
and dinaturality.
\end{proof}

The following theorem now extends Quillen's criterion to
$\lambda$-enriched functors. The corresponding criterion for right
derived functors also holds (and follows by considering
$\Phi^{\op}\colon \EC^{\op}\to\EH^{\op}$).

\begin{theorem}[Enriched Quillen Criterion]\label{enrquilcri}
Let $\EC$ be an enriched model category over a monoidal model category
$\Moncat$, $\EH$ be a category enriched over $\Ho\Moncat$, and $\Phi
\colon \EC\to \EH$ be a $\Ho\Moncat$-enriched functor.  If $\Phi$
takes weak equivalences between cofibrant objects to isomorphisms in
$\EH$, then the enriched left derived functor exists and provides an
enrichment for the left derived functor.
\end{theorem}

\begin{proof}
The composite enriched functor $\Phi \circ Q\colon \lambda_{*}\EC\to
\EH$ sends weak equivalences to isomorphisms, and so factors through
an enriched functor $L^{\Moncat}\Phi \colon \Ho\EC\to \EH$, and
$\phi^{\Moncat} =\Phi(q)$ gives a natural transformation from
$L^{\Moncat}\Phi \circ \gamma=\Phi \circ Q$ to $\Phi$.  It is easy to
see from the diagrams preceding Lemma~\ref{Qenr} that the underlying
functor and natural transformation are the left derived functor and
universal natural transformation constructed by Quillen.  Given any
enriched functor $F\colon \Ho\EC\to \EH$ and any enriched natural
transformation $\alpha\colon F\circ \gamma \to \Phi$, the maps
$\theta_{\objectone} = \alpha_{Q\objectone} \circ
F(q^{-1}_{\objectone})$ assemble to a natural transformation $\theta
\colon F\to L\Phi$, and this is the unique natural transformation
$\theta$ such that $\phi^{\Moncat} \circ \theta =\alpha$.
\end{proof}


Next we discuss Quillen adjunctions.  Recall that given an adjunction
between closed model categories $\EC$ and $\ED$, the following are
equivalent:
\begin{hyplist}
\item The left adjoint preserves cofibrations and acyclic
cofibrations.
\item The right adjoint preserves fibrations and acyclic fibrations.
\item The left adjoint preserves cofibrations and the right adjoint
preserves fibrations.
\end{hyplist}
Such an adjunction is called a \defterm{Quillen adjunction}.  The left adjoint
$\Phi \colon \EC\to \ED$ then preserves weak equivalences between
cofibrant objects, and the right adjoint $\Theta \colon \ED\to \EC$
preserves weak equivalences between fibrant objects.  Quillen's
criterion for the existence of left derived functors then applies to
$\gamma_{\ED}\circ \Phi$ and Quillen's criterion for the existence of
right derived functors to $\gamma_{\EC}\circ \Theta$.  We write
$\tL\Phi$ and $\tR\Theta$ for the corresponding (unenriched and, if
applicable, enriched) derived functors obtained from
$\gamma_{\ED}\circ \Phi$ and $\gamma_{\EC}\circ \Theta$.  The
fundamental theorem of model category theory is that $\tL\Phi$ and
$\tR\Theta$ remain adjoints.  To extend this to the enriched context
we first must recall the definition of an enriched adjunction.

\begin{definition}
An adjunction between enriched functors $\Phi\colon \EC \to
\ED$ and $\Theta\colon \ED \to \EC$ is said to be \defterm{enriched}
if the unit and counit of the adjunction are both enriched natural.
This condition is equivalent to the requirement that the adjunction
isomorphism
\[
\EC(C,\Theta D) \iso \ED(\Phi C,D)
\]
lifts to an isomorphism
\[
\efunc[{\EC}]{C}{\Theta D} \iso \efunc[{\ED}]{\Phi C}{D}.
\]
that is enriched natural in each variable.
\end{definition}

\begin{theorem}[Enriched Quillen Adjunction]\label{enradj}
Let $\EC$ and $\ED$ be enriched model categories over the monoidal
model category $\Moncat$.  If an enriched adjunction $(\Phi,\Theta)$
between $\EC$ and $\ED$ is a Quillen adjunction, then the derived
adjunction $(\tL\Phi,\tR\Theta)$ is also enriched.
\end{theorem}

\begin{proof}
If we write $\eta \colon \Id\to \Theta \Phi$ for the unit of the
$(\Phi,\Theta)$ adjunction, then the unit of the $(\tL\Phi
,\tR\Theta)$ adjunction is $(\Theta r)\circ \eta \circ q^{-1}$, and
this is clearly enriched when $\eta$ is, and likewise for the counit.
\end{proof}

As promised in the last section, we now complete the proof of
Theorem~\ref{mainenr}.  This amounts to recalling the notions of
tensors and cotensors and applying the result above.

\begin{definition}\label{TcTDef}
For an object $\objectone$ of $\EC$ and an object $X$ of $\Moncat$,
the associated \defterm{tensor} $\tensor{\objectone}{X}$ and \defterm{cotensor}
$\cotensor{X}{\objectone}$ are objects of $\EC$, unique up to an
enriched natural isomorphism when they exist, for which there are
enriched natural isomorphisms
\[
\efunc{\tensor{\objectone}{X}}{-}\iso \func{X}{\efunc{\objectone}{-}}
\tand \efunc{-}{\cotensor{X}{\objectone}}\iso
\func{X}{\efunc{-}{\objectone}}.
\]
If $\tensor{\objectone}{X}$ exists for all $X$, then for formal
reasons $\tensor{\objectone}{(-)}$ is an enriched functor, and we can
interpret the natural isomorphism above as an enriched adjunction.
Analogous observations hold for $\cotensor{-}{\objectone}$.
\end{definition}

Applying the previous theorem to these enriched adjunctions gives the
following corollary and thereby completes the proof of
Theorem~\ref{mainenr}.

\begin{corollary}\label{mainenrtwo}
Let $\EC$ be an enriched model category over a monoidal model category
$\Moncat$.
\begin{thmlist}
\item If $\EC$ has tensors then so does $\Ho\EC$ and the tensor in
$\Ho\EC$ is the left derived functor $\emptydertensor$ of the tensor
$\emptytensor$ in $\EC$
\item If $\EC$ has cotensors then so does $\Ho\EC$ and the cotensor in
$\Ho\EC$ is the right derived functor $\emptydercotensor$ of the
cotensor $\emptycotensor$ in $\EC$.
\end{thmlist}
\end{corollary}


\section{Enriched bifunctors and their derived
functors}\label{BiFunSec}

In the context of ordinary category theory, bifunctors such as the
functor taking a pair of objects $C$ and $D$ in a category $\EC$ to
their product $C \times D$ or to the morphism set $\EC(C,D)$ have as
their domains categories of the form $\EC \times \EC$, $\EC \times
\EC^\op$, or $\EC \times \ED$.  However, the domains of the analogous
enriched bifunctors, such as the tensor product functor, enriched hom
functors, and the tensor and cotensor functors are not product
categories like $\EC \times \ED$, but more complex enriched categories
of the form $\EC \monprod \ED$, as can be seen in the familiar
examples of additive categories.  In the context of enriched model
categories, this problem with domain categories for bifunctors is
compounded by the fact that the morphism sets of the ordinary category
underlying an enriched category like $\EC \monprod \ED$ typically have
no tractable description.  As a result, we cannot expect to be able to
impose a useful model structure on these categories.  The purpose of
this section is to propose a definition of enriched derived functors
in this context and to study when they exist and fit into
(parametrized) enriched adjunctions.

We begin by reviewing the definition of enriched bifunctor.  For
categories $\EC$ and $\ED$ enriched over $\Moncat$, the enriched
category $\EC \monprod \ED$ is defined to have objects
$\ob{\EC\monprod\ED}=\ob{\EC}\times \ob{\ED}$ and for objects
$(C,D),(C',D')$ in $\ob{\EC\monprod\ED}$, the morphism object
$\efunc[\EC\monprod\ED]{(C,D)}{(C',D')}$ in \Moncat is defined to be
\begin{equation*}
\efunc[\EC]{C}{C'}\monprod \efunc[\ED]{D}{D'}.
\end{equation*}
An enriched bifunctor from $\EC,\ED$ to an enriched category $\EE$ is
defined to be an enriched functor $\EC\monprod \ED\to \EE$.  The
ordinary bifunctor $\EC \times \ED \to \EE$ underlying an enriched
bifunctor $\EC \monprod \ED \to \EE$ is obtained by precomposing with
the functor from from $\EC\times \ED$ to the underlying category of
$\EC\monprod \ED$ that takes
\[
f\in \Moncat(\unit,\efunc[\EC]{C}{C'})=\EC(C,C'), \qquad g\in
\Moncat(\unit,\efunc[\ED]{D}{D'})=\ED(D,D')
\]
to
\[
f\monprod g \in \Moncat(\unit,\efunc[\EC]{C}{C'} \monprod
\efunc[\ED]{D}{D'}) = (\EC\monprod \ED)((C,D),(C',D')).
\]
An enriched natural transformation of bifunctors is an enriched
natural transformation of functors $\EC\monprod \ED\to \EE$, or
equivalently, a natural transformation that is enriched in each
variable separately.

In the context of a monoidal model category \Moncat, since $\lambda
\colon \Moncat\to \Ho\Moncat$ is lax symmetric monoidal, we have a
canonical $\Ho\Moncat$-enriched functor
\[
\dermonprod{\lambda_{*}\EC}{\lambda_{*}\ED}\to
\lambda_{*}(\EC\monprod\ED).
\]
This functor is typically not an equivalence.  When $\EC$ and $\ED$
are enriched model categories, the $\Ho\Moncat$-enriched localization
functors $\gamma_{\EC} \colon \lambda_{*}\EC\to \Ho\EC$ and
$\gamma_{\ED} \colon \lambda_{*}\ED\to \Ho\ED$ induce a
$\Ho\Moncat$-enriched functor
\[
\dermonprod{\gamma_{\EC}}{\gamma_{\ED}} \colon
\dermonprod{\lambda_{*}\EC}{\lambda_{*}\ED} \to
\dermonprod{\Ho\EC}{\Ho\ED},
\]
but we do not expect an enriched functor from
$\lambda_{*}(\EC\monprod\ED)$ to $\dermonprod{\Ho\EC}{\Ho\ED}$.  The
following theorem describing the universal property of
$\dermonprod{\gamma_{\EC}}{\gamma_{\ED}}$ is the bifunctor analog of
Theorem \ref{localenr}; its proof is a straightforward application of
diagram~\eqref{enrcompdiag}.

\begin{theorem}\label{bilocalenr}
Let $\EC$ and $\ED$ be enriched model categories over the monoidal
model category $\Moncat$.  The $\Ho\Moncat$-enriched bifunctor
\map{\dermonprod{\gamma_{\EC}}{\gamma_{\ED}}}
{\dermonprod{\lambda_*\EC}{\lambda_*\ED}}{\dermonprod{\Ho\EC}{\Ho\ED}}
is initial among $\Ho\Moncat$-enriched bifunctors that take weak
equivalences in each variable to isomorphisms; in other words, if
$\EH$ is a $\Ho\Moncat$-enriched category and \map{\Phi}{
\dermonprod{\lambda_*\EC}{\lambda_*\ED}}{\EH} is a
$\Ho\Moncat$-enriched functor that sends weak equivalences to
isomorphisms, then $\Phi$ factors uniquely through a
$\Ho\Moncat$-enriched functor
\[ \vmap{\dermonprod{\Ho\EC}{\Ho\ED}}{\EH}. \]
\end{theorem}

Next we discuss the enriched left derived functors of enriched
bifunctors.  We offer the following definition.

\begin{definition}
Let $\EC$ and $\ED$ be enriched model categories over a monoidal model
category $\Moncat$, let $\EH$ be a category enriched over $\Ho\Moncat$,
and let $\Phi \colon \dermonprod{\lambda_*\EC}{\lambda_*\ED} \to \EH$ be a
$\Ho\Moncat$-enriched bifunctor.  An enriched bifunctor
\map{L^{\Moncat}\Phi}{\dermonprod{\Ho\EC}{\Ho\ED}}{\EH} with an
enriched natural transformation $\phi^{\Moncat} \colon L^{\Moncat}\Phi
\circ (\dermonprod{\gamma_{\EC}}{\gamma_{\ED}}) \to \Phi$ forms the
\defterm{enriched left derived bifunctor} of $\Phi$ when $\phi^{\Moncat}$ is
final among enriched natural transformations $F\circ
(\dermonprod{\gamma_{\EC}}{\gamma_{\ED}})\to\Phi$.  We say that the
enriched left derived bifunctor enriches the left derived functor of
$\Phi$ if the left derived functor $L\Phi,\phi$ of $\Phi$ exists and
is the restriction to $\Ho\EC\times \Ho\ED$ of the underlying functor
and natural transformation of $L^{\Moncat}\Phi,\phi^{\Moncat}$.
\end{definition}

Enriched right derived bifunctors are defined analogously.  The
following theorem is the bifunctor equivalent of Theorem
\ref{enrquilcri}.  The corresponding result for right derived functors
also holds (and follows by considering the appropriate enriched
opposite categories).

\begin{theorem}\label{penrquilcri}
Let $\EC$ and $\ED$ be enriched model categories over a monoidal model
category $\Moncat$, let $\EH$ be a category enriched over $\Ho\Moncat$,
and let $\Phi \colon \dermonprod{\lambda_*\EC}{\lambda_*\ED} \to \EH$ be a
$\Ho\Moncat$-enriched bifunctor.  If $\Phi$ takes weak equivalences
between cofibrant objects to isomorphisms in $\EH$, then the enriched
left derived bifunctor exists and enriches the left derived functor.
\end{theorem}

\begin{proof}
We obtain $L^{\Moncat}\Phi$ by factoring $\Phi \circ
(\dermonprod{Q}{Q})$ using Theorem~\ref{bilocalenr}.  Given an
enriched functor $F\colon \dermonprod{\Ho\EC}{\Ho\ED}\to \EH$ and an
enriched natural transformation $\alpha \colon F\circ
\dermonprod{\gamma_{\EC}}{\gamma_{\ED}}\to \Phi$, then
$F(\dermonprod{q^{-1}}{q^{-1}})$ is the unique enriched natural
transformation factoring $\alpha$.
\end{proof}

The definition and theorems above have obvious generalizations to
trifunctors and functors of any number of variables.  In general for
$\EC_{0},\ldots,\EC_{m}$ enriched model categories over a monoidal
model category $\Moncat$, and an enriched functor of $m$-variables
\[
\Phi \colon \EC_{1}\monprod \cdots \monprod \EC_{m}\to \EC_{0},
\]
we write
\[
\tL\Phi \colon \Ho\EC_{1} \bindmp \cdots \bindmp \Ho\EC_{m}\to
\Ho\EC_{0}
\]
for the enriched left derived functor of
\[
\gamma_{\EC_{0}}\circ \lambda_{*}\Phi \colon \lambda_{*}\EC_{1}
\bindmp \cdots \bindmp \lambda_{*}\EC_{m}\to \lambda_{*}\EC_{0} \to
\Ho\EC_{0}
\]
when it exists and extends the left derived functor, and call it the
\defterm{enriched total left derived functor}.  The enriched total right derived
functor is defined analogously and denoted $\tR\Phi$.  The following
terminology is also convenient.

\begin{definition}\label{totalderivedfunctor}
For $\Phi$ as above, we say that the left derived functor of $\Phi$ is
\defterm{enriched} when the enriched left derived functor of $\gamma\circ
\lambda_{*}\Phi$ exists and extends the left derived functor.
Likewise, we say that the right derived functor of $\Phi$ is enriched
when the enriched right derived functor of $\gamma\circ
\lambda_{*}\Phi$ exists and extends the right derived functor.
\end{definition}

Functors of many variables admit many sorts of compositions, and the
same kind of results as usual for the composition of derived functors
of a single variable apply to all the possible compositions of total
derived functors of many variables.  The following proposition
suffices for our purposes in Section~\ref{PrfSec}.  We phrase the
proposition for the enriched total derived functors but it is really
an assertion about the unenriched total derived functors.

\begin{proposition}\label{composederived}
Let $\EC_{0},\ldots,\EC_{m}$, $\ED_{0},\ldots,\ED_{n},$ and $\EE$ be
enriched model categories over a monoidal model category $\Moncat$. If
\[
\Phi \colon \EC_{0} \monprod \ED_{0}\to\EE,\qquad \Psi \colon \EC_{1}
\monprod \cdots \monprod \EC_{m}\to \EC_{0},\qquad \Xi \colon \ED_{1}
\monprod \cdots \monprod \ED_{n}\to \ED_{0},\qquad
\]
are enriched functors that send tuples of cofibrant objects to
cofibrant objects and preserve weak equivalences between tuples of
cofibrant objects, then the universal map
\[
(\tL\Phi) \circ (\dermonprod{\tL\Psi}{\tL\Xi}) \to \tL(\Phi \circ
({\Psi }\monprod{\Xi}))
\]
is an isomorphism.
\end{proposition}


All of the enriched bifunctors of interest to us appear in enriched
parametrized adjunctions, and it is important that these adjunctions
pass to homotopy categories.  The general context we study is when
have a pair of bifunctors
\begin{equation*}
\Phi \colon \EC \monprod \ED \to \EE \tand \Theta \colon
\ED^{\op}\monprod \EE\to \EC
\end{equation*}
that form an enriched parametrized adjunction.  This means that we
have isomorphisms
\begin{equation*}
\efunc[\EE]{\Phi(C,D)}{E}\iso \efunc[\EC]{C}{\Theta(D,E)}.
\end{equation*}
that are enriched natural in all three variables.  The following
proposition describes the two pairs of equivalent conditions that
together suffice to ensure that such a parametrized adjunction passes
properly to homotopy categories.

\begin{proposition}\label{ParmAdjImp}
Let $\EC$, $\ED$, and $\EE$ be enriched model categories over a
monoidal model category \Moncat, and let $\Phi \colon \EC \monprod \ED
\to \EE$ and $\Theta \colon \ED^{\op}\monprod \EE\to \EC$ be a pair of
\Moncat-enriched bifunctors forming an enriched parametrized adjoint
pair.  The following two conditions are equivalent:
\begin{hyplist}
\item $\Phi(-,D)$ preserves cofibrations and acyclic cofibrations for
all cofibrant $D$ in $\ED$

\item $\Theta(D,-)$ preserves fibrations and acyclic fibrations for
all cofibrant $D$ in $\ED$
\end{hyplist}
If $\Phi$ and $\Theta$ satisfy these conditions, then the following
two conditions are equivalent:
\begin{thmlist}
\item $\Phi (C,-)$ preserves weak equivalences between cofibrant
objects for all cofibrant $C$ in $\EC$.

\item $\Theta(-,E)$ preserves weak equivalences between cofibrant
objects for all fibrant $E$ in $\EE$
\end{thmlist}

\end{proposition}

\begin{proof}
The equivalence of (i) and (ii) is just a special case of the standard
result about Quillen adjunctions.  Assuming (i) and (ii), let
\map{f}{D}{D'} be a weak equivalence between two cofibrant objects in
$\ED$, $C$ be a cofibrant object of $\EC$, and $E$ be a fibrant object
of $\EE$.  Then $\Phi(C,D)$ and $\Phi(C,D')$ are cofibrant and
$\Theta(D,E)$ and $\Theta(D',E)$ are fibrant by (i) and (ii).  It
follows that we can identify the commuting diagram on the left below
with the commuting diagram on the right below.
\begin{equation*}
\xymatrix@C-.5ex{%
{\efunc[\EE]{\Phi(C,D')}{E} \ar[r]^-{\iso}
\ar[d]_{\Phi(\id_{C},f)^*} } & {\efunc[\EC]{C}{\Theta(D',E)}
\ar[d]^{\Theta(f,\id_{E})_*}} &
{\derefunc[\EE]{\Phi(C,D')}{E} \ar[r]^-{\iso} \ar[d]_{\Phi(\id_{C},f)^*} } & {\derefunc[\EC]{C}{\Theta(D',E)} \ar[d]^{\Theta(f,\id_{E})_*}}  \\
{\efunc[\EE]{\Phi(C,D)}{E} \ar[r]^-{\iso}} &
{\efunc[\EC]{C}{\Theta(D,E)}}& {\derefunc[\EE]{\Phi(C,D)}{E}
\ar[r]^-{\iso}} & {\derefunc[\EC]{C}{\Theta(D,E)}} }
\end{equation*}
Then $\Phi(\id_{C},f)^{*}$ is an isomorphism for every
cofibrant $C$ in $\EC$ and every fibrant $E$ in $\EE$ if and only if
$\Theta(f,\id_{E})_{*}$ is an isomorphism for every cofibrant $C$ in
$\EC$ and every fibrant $E$ in $\EE$.  Now by the enriched Yoneda
Lemma in $\Ho\EC$ and $\Ho\EE$, we see that $\Phi(\id_{C},f)$ is a
weak equivalence for every cofibrant $C$ if and only if
$\Theta(f,\id_{E})$ is a weak equivalence for every fibrant $E$ in
$\EE$.
\end{proof}


We can now state our main result on the passage of parametrized
adjunctions to homotopy categories.

\begin{theorem}\label{epqa}
Let $(\Phi,\Theta)$ be an enriched parametrized adjunction satisfying
both pairs of equivalent conditions of
Proposition~\ref{ParmAdjImp}. Then the enriched total left derived
bifunctor
\begin{equation*}
\map{\tL\Phi}{\dermonprod{\Ho\EC}{\Ho\ED}}{\Ho\EE}
\end{equation*}
of $\Phi$ and the enriched total right derived bifunctor
\begin{equation*}
\map{\tR\Theta}{\dermonprod{\Ho\ED^{\op}}{\Ho\EE}}{\Ho\EC}
\end{equation*}
of $\Theta$ exist and form an enriched parametrized adjunction on the
homotopy categories.
\end{theorem}

\begin{proof}
The enriched total derived bifunctors exist by
Theorem~\ref{penrquilcri} and its analogue for right derived functors.
For fixed $D$ in $\ED$, write $\Phi_{D}$ for the enriched functor
$\Phi(-,QD)\colon \EC\to \EE$ and $\Theta_{D}$ for the enriched
functor $\Theta(QD,-)\colon \EE\to \EC$.  Then by
Theorem~\ref{enradj}, the total derived functors $\tL\Phi_{D}\colon
\Ho\EC\to \Ho\EE$ and $\tR\Theta_{D}\colon \Ho\EE\to \Ho\EC$ are
adjoint.  By construction, $\tL\Phi_{D}$ coincides with $\tL\Phi
(-,D)$ viewed as an enriched functor $\Ho\EC\to\Ho\EE$, and
$\tR\Theta_{D}$ coincides with $\tR\Theta(D,-)$ viewed as an enriched
functor $\Ho\EE\to\Ho\EC$.  Viewing $\tL\Phi_{(-)}(-)$ as the
bifunctor $\tL\Phi \colon \dermonprod{\Ho\EC}{\Ho\ED}\to \Ho\ED$, then
formally there exists precisely one way to make $\tR\Theta_{(-)}$ an
enriched bifunctor $\dermonprod{\Ho\ED^{\op}}{\Ho\EE}\to \Ho\EC$ that
is a parametrized right adjoint to $\tL\Phi$; we have to show that for
fixed $E$ in $\EE$, the functor $\tR\Theta_{(-)}(E)\colon
\Ho\ED^{\op}\to \Ho\EC$ so obtained coincides with $\tR\Theta(-,E)$.
Denote the counit of the $\tL\Phi_{D'}, \tR\Theta_{D'}$ adjunction as
$\epsilon'$ and the natural isomorphism
$\derefunc{\tL\Phi_{D}(-)}{-}\iso \derefunc[\EE]{-}{\tR\Theta_{D}(-)}$
as $\alpha'$.  Then as a contravariant $\Ho\Moncat$-enriched functor
on $\Ho\ED$, $\tR\Theta_{(-)}(E)$ is the map
\begin{multline*}
\derefunc[\ED]{D}{D'}\lrarrow[\tL\Phi(\tR \Theta_{D'}(E),-)]
\derefunc{\tL\Phi(\tR\Theta_{D'}(E),D)}{\tL\Phi(\tR\Theta_{D'}(E),D')}\\
\lrarrow[\epsilon'_{*}] \derefunc{\tL\Phi(\tR\Theta_{D'}(E),D)}{E}
\lrarrow[\alpha'] \derefunc[\EE]{\tR\Theta_{D'}(E)}{\tR\Theta_{D}(E)}
\end{multline*}
Unwinding the definition of $\tL\Phi$, $\tR\Theta$, and using
naturality in $\ED$ of the $\Phi,\Theta$ adjunction, a little bit of
work identifies this map as the composite in $\Ho\Moncat$
\begin{multline*}
\derefunc[\ED]{D}{D'}= \efunc[\ED]{QD}{RD'}\simeq\efunc[\ED]{QD}{RQD}
\lrarrow[\Phi(Q\Theta(RQD',RE)),Q(-))]\\
\efunc{\Phi(Q\Theta(RQD',RE),QD)}{\Phi(Q\Theta(RQD',RE),RQD')}
\lrarrow[\epsilon_{*}]\qquad\\
\qquad\efunc{\Phi(Q\Theta(RQD',RE),QD)}{RE} \lrarrow[\alpha]
\efunc[\EE]{Q\Theta(RQD',RE)}{\Theta(QD,RE)}\\
\simeq \efunc[\EE]{Q\Theta(QD',RE)}{\Theta(QD,RE)}
=\derefunc[\EE]{\tR\Theta(D',E)}{\tR\Theta(D,E)},
\end{multline*}
where $\epsilon$ is the counit and $\alpha$ the isomorphism for the
$\Phi,\Theta$ adjunction.  Unwinding the $\Phi,\Theta$ adjunction
identifies this composite as the functor $\tR\Theta$.
\end{proof}


\section{Semicofibrant objects}\label{SecSemi}

As explained in the introduction, semicofibrant objects in $\Moncat$
are of intrinsic interest because in practice monoids in $\Moncat$ can
often be approximated by weakly equivalent monoids whose underlying
objects are semicofibrant, but when the unit $\unit$ is not cofibrant,
monoids typically cannot have underlying objects that are cofibrant.
We need some further observations on the properties of semicofibrant
objects for the proofs of the main results of the introduction that
are phrased in terms of semicofibrant objects.  In this section, we
collect these observations and some additional facts about
semicofibrant objects that seem potentially useful.

Recall that an object $C$ in a closed model category $\EC$ enriched
over the monoidal model category \Moncat is semicofibrant when the
functor $\efunc{C}{-}\colon \EC\to \Moncat$ preserves fibrations and
acyclic fibrations.  Clearly, this notion is most useful when $\EC$ is
an enriched model category, and we have the following proposition that
generalizes parts of Proposition~\ref{SemiCoPropA}.

\begin{proposition}\label{semigen}
Let $\EC$ be an enriched model category.
\begin{thmlist} 

\item If $C$ is cofibrant in \EC, then it is semicofibrant in
\EC. Moreover, if the unit \unit is cofibrant in \Moncat, then an
object $C$ of \EC is semicofibrant in \EC if and only if it is
cofibrant in \EC.

\item If $C$ is semicofibrant in \EC and $C\to D$ is a cofibration in
\EC, then $D$ is semicofibrant in \ED.
\end{thmlist}
\end{proposition}

\begin{proof}
The first part of part~(a) is a special case of part~(b), which follows
immediately from the Enrichment Axiom.  For the second part of
part~(a), suppose \unit is cofibrant in \Moncat and let $C$ be a
semicofibrant object in $\EC$; then to see that $C$ is cofibrant, we
just need to see that for any acyclic fibration $X\to Y$ and any map
$C\to Y$, there exists a lift $C\to X$.  A map $C\to Y$ specifies a
map in \Moncat from \unit to $\efunc{C}{Y}$.  Since $C$ is
semicofibrant, $\efunc{C}{X}\to \efunc{C}{Y}$ is an acyclic fibration
in \Moncat.  Since \unit is cofibrant in \Moncat, we can lift
$\unit\to \efunc{C}{Y}$ to $\unit\to\efunc{C}{X}$, and this specifies
the lift $C\to X$ in $\EC$ of $C\to Y$.
\end{proof}

We also need the following general theorem about semicofibrant
objects.  It is proved at the end of the section.

\begin{theorem}\label{semigood}
Let $\EC$ be an enriched model category and $D$ a fibrant object of
$\EC$.  Then $\efunc{-}{D}$ preserves weak equivalences between
semicofibrant objects.
\end{theorem}

Applying the theorem to the cofibrant approximation $QC\to C$, we
obtain the following corollary.

\begin{corollary}\label{corsemigood}
Under the hypotheses of Theorem~\ref{semigood}, if $C$ is a
semicofibrant object of $\EC$, then the canonical map $\efunc{C}{D}\to
\derefunc{C}{D}$ is an isomorphism in $\Ho\Moncat$.
\end{corollary}

The following proposition explains many of the properties of
semicofibrant objects.

\begin{proposition}\label{semicosmash}
Let $\EC$ be an enriched model category that has tensors.
\begin{thmlist}
\item $C$ is semicofibrant if and only if $\tensor{C}{-}\colon
\Moncat\to \EC$ preserves cofibrations and acyclic cofibrations.
\item If $C$ is semicofibrant, then $\tensor{C}{\cofunit}$ is
cofibrant, and $\cum\colon \tensor{C}{\cofunit}\to C$ is a weak
equivalence.
\end{thmlist}
\end{proposition}

\begin{proof}
Part~(a) is the usual usual result on Quillen adjunctions applied to
the adjoint pair $\tensor{C}{-}\colon \Moncat\to \EC$ and
$\efunc{C}{-}\colon \EC\to \Moncat$.  The first statement of part~(b)
follows from applying part~(a) to the cofibration $0\to \cofunit$
where $0$ is the initial object.  For the second statement in
part~(b), consider the map
\[
\efunc{C}{D}\to \func{\cofunit}{\efunc{C}{D}}\iso
\efunc{\tensor{C}{\cofunit}}{D}.
\]
When $D$ is fibrant in $\EC$, $\efunc{C}{D}$ is fibrant in $\Moncat$,
and so this map is a weak equivalence by the Unit Axiom in $\Moncat$.
Since the composite is induced by the map $\cum\colon
\tensor{C}{\cofunit}\to C$, applying Corollary~\ref{corsemigood} and
the enriched Yoneda Lemma, we see that $\cum$ is a weak equivalence.
\end{proof}

Finally, we need the following two propositions which are specific to
the case of module categories.  The first proposition is clear from
the definition of semicofibrant. Together with it,
Propositions~\ref{semigen} and~\ref{semicosmash} subsume
Proposition~\ref{SemiCoPropA} from the introduction.

\begin{proposition}
Assume that \Leftmod is a closed model category with fibrations and
weak equivalences created in \Moncat.  Then \monoid, considered as an
object of \Leftmod, is semicofibrant in \Leftmod.
\end{proposition}

The following proposition is a formal statement of the observation
following Theorem~\ref{intexist} and plays a key role in the arguments
in the next section.

\begin{proposition}\label{SemiCoPropB} 
Assume that \Leftmod and \Bimod{\monoid}{\monoidtwo} are closed model
categories with fibrations and weak equivalences created in \Moncat.
If the monoid $\monoidtwo$ is semicofibrant as an object of $\Moncat$,
then the forgetful functor $\Bimod{\monoid}{\monoidtwo}\to \Leftmod$
preserves cofibrations and takes semicofibrant objects in
$\Bimod{\monoid}{\monoidtwo}$ to semicofibrant objects in $\Leftmod$.
\end{proposition}

\begin{proof}
The functor
\map{\func{\monoidtwo}{-}}{\Leftmod}{\Bimod{\monoid}{\monoidtwo}} is
right adjoint to the forgetful functor from
\Bimod{\monoid}{\monoidtwo} to \Leftmod.  Since $\monoidtwo$ is
semicofibrant in \Moncat, this right adjoint preserves fibrations and
acyclic fibrations.  Thus, the forgetful functor from
\Bimod{\monoid}{\monoidtwo} to \Leftmod preserves cofibrations (and by
hypothesis on the model structures, all weak equivalences).  If $M$ is
semicofibrant in \Bimod{\monoidone}{\monoidtwo}, then we see from the
natural isomorphism
\[
\lmfunc{M}{-}\iso \bmfunc{\monoid}{\monoidtwo}{M}{\func{B}{-}}\colon
\Leftmod\to\Moncat
\]
that $\lmfunc{M}{-}$ preserves fibrations and acyclic fibrations, and
so $M$ is semicofibrant in \Leftmod.
\end{proof}

We close the section with the proof of Theorem~\ref{semigood}.  Let
$\EC$ be an enriched model category, let $D$ be a fibrant object, and
let $f\colon C'\to C$ be a weak equivalence between semicofibrant
objects; we need to see that $f^{*}\colon \efunc{C}{D}\to
\efunc{C'}{D}$ is a weak equivalence.  Since $\efunc{-}{D}$ converts
acyclic cofibrations to acyclic fibrations, by factoring the map from
$C$ to the final object, it suffices to consider the case when $C$ is
fibrant.  Likewise, by factoring the map $f$, it suffices to consider
the case when $f$ is an acyclic fibration.

The idea for the proof is to construct some kind of ``map'' $g\colon
C\to C'$ such that the composite $f\circ g\colon C\to C$ is the
identity and the composite $g\circ f\colon C'\to C'$ is (left)
homotopic to the identity.  The induced composite $g^{*}\circ f^{*}$
then would be the identity and $f^{*}\circ g^{*}$ would be (right)
homotopic to the identity, and so still a weak equivalence.  We can
actually do this in the case when \unit is cofibrant using a version
of the argument of Proposition~\ref{semigen} (or the proposition
itself).  We generalize this argument and make this idea rigorous as
follows:

Since $C$ is semicofibrant and $f$ is an acyclic fibration, the map
$f_{*}\colon \efunc{C}{C'}\to \efunc{C}{C}$ is an acyclic fibration.
Let $\tilde{\id}_{C} \colon \cofunit \to \efunc{C}{C}$ be the
composite of $\cofunit \to \unit$ and the map $\unit\to\efunc{C}{C}$
representing the identity of $C$.  Then since $\cofunit$ is cofibrant,
we can lift $\tilde{\id}_{C}$ to a map $g\colon \cofunit \to
\efunc{C}{C'}$.
\[
\xymatrix{%
\relax&\relax\efunc{C}{C'}\ar@{->>}[d]_{\sim}^{f_{*}}\\
\cofunit \ar[r]_{\tilde{\id}_{C}}\ar@{-->}[ur]^{g}&\relax\efunc{C}{C}
}
\]
Composition gives a map
\[
\efunc{C'}{D} \monprod \cofunit \to \efunc{C'}{D} \monprod
\efunc{C}{C'} \to \efunc{C}{D}
\]
and adjoint to this map, we have a map
\[
\hat{g}\colon \efunc{C'}{D}\to \func{\cofunit}{\efunc{C}{D}}.
\]
By construction, the composite map $\hat{g}\circ f^{*}\colon
\efunc{C}{D}\to \func{\cofunit}{\efunc{C}{D}}$ is the map $\cuf$,
which is a weak equivalence by the Unit Axiom in \Moncat (since
$\efunc{C}{D}$ is fibrant).

We have constructed a commutative diagram
\[
\xymatrix{%
\efunc{C}{D}\ar[r]^{\cuf}_{\sim}\ar[d]_{f^{*}}
&\func{\cofunit}{\efunc{C}{D}}\ar[d]^{(f^{*})_{*}}\\
\efunc{C'}{D}\ar[r]\ar[ur]^{\hat{g}}
&\func{\cofunit}{\efunc{C'}{D}}
}
\]
where the bottom map is the composite $(f^{*})_{*}\circ \hat{g}$.  If
this map were $\cuf$, then the argument would be complete
(cf. Proposition~\ref{proptwosix} below); the remainder of the
argument is to show that it is homotopic to $\cuf$.  Note that
$(f^{*})_{*}\circ \hat{g}$ is 
induced by ``composition'' with the map $h=f^{*}\circ g\colon \cofunit
\to \efunc{C'}{C'}$.  If we write $\tilde{f}$ for the map $\cofunit
\to \efunc{C'}{C}$ adjoint to the map $\cofunit\to \unit \to
\efunc{C'}{C}$ representing $f$, then we have
\[
f_{*}\circ h=f_{*}\circ f^{*}\circ g = f^{*}\circ \tilde{\id}_{C}
=\tilde{f}\colon \cofunit \to \efunc{C'}{C}.
\]
Let $J$ be a Quillen left cylinder object for $\cofunit$, i.e., factor
the codiagonal map $\cofunit \amalg \cofunit\to \cofunit$ as a
cofibration $\cofunit \amalg \cofunit\to J$ followed by an acyclic
fibration $J\to \cofunit$.  Then writing $\tilde{\id}_{C'}\colon
\cofunit\to \efunc{C'}{C'}$ for the map $\cofunit \to \unit\to
\efunc{C'}{C'}$ induced by the identity on $C'$, we have the following
solid arrow commuting diagram:
\[
\xymatrix{%
\relax\cofunit \amalg \cofunit
\ar[rr]^{\widetilde{\id}_{C'}\amalg h}\ar@{>->}[d]
&&\relax\efunc{C'}{C'}\ar@{->>}[d]_{\sim}^{f_{*}}\\
J\ar@{-->}[urr]^{\phi}\ar[r]&\cofunit\ar[r]_(.4){\tilde{f}}
&\relax\efunc{C'}{C} }
\]
Choose a lift $\phi$ as indicated by the dashed arrow in the diagram.
Then composition gives us a map
\[
\efunc{C'}{D}\monprod J\to \efunc{C'}{D}\monprod \efunc{C'}{C'}\to
\efunc{C'}{D}.
\]
Let $\hat{\phi}$ denote the adjoint map $\efunc{C'}{D}\to
\func{J}{\efunc{C'}{D}}$.

The two acyclic cofibrations $\cofunit \to J$ induces two acyclic
fibrations $\func{J}{\efunc{C'}{D}}\to
\func{\cofunit}{\efunc{C'}{D}}$.  By composition with $\hat{\phi}$, we
obtain two maps $\efunc{C'}{D}\to \func{\cofunit}{\efunc{C'}{D}}$,
which by construction are $\cuf$ and $(f^{*})_{*}\circ \hat{g}$.
It follows that $\func{J}{\efunc{C'}{D}}$ is a Quillen path object for
$\func{\cofunit}{\efunc{C'}{D}}$ and that $\hat{\phi}$ is a Quillen right
homotopy between $(f^{*})_{*}\circ \hat{g}$ and $\cuf$.

In particular, since $\cuf$ is a weak equivalence,  $\hat{\phi}$ is a
weak equivalence, and therefore $(f^{*})_{*}\circ \hat{g}$ is a
weak equivalence.
This shows that of the three composable maps
\begin{align*}
f^{*}&\colon \efunc{C}{D}\to \efunc{C'}{D}\\
\hat{g}&\colon \efunc{C'}{D}\to \func{\cofunit}{\efunc{C}{D}}\\
(f^{*})_{*}&\colon \func{\cofunit}{\efunc{C}{D}}\to
\func{\cofunit}{\efunc{C'}{D}}
\end{align*}
both $\hat{g}\circ f^{*}$ and $(f^{*})_{*}\circ \hat{g}$ are weak
equivalences.  The following ``two out of six'' principle of Dwyer,
Hirschhorn, Kan, and Smith \cite[8.2.(ii)]{DHKS} implies that $f^{*}$
is a weak equivalence and completes the proof of the theorem.

\begin{proposition}[Two out of Six Principle {\cite[\S9]{DHKS}}]\label{proptwosix}
Let $\EC$ be a closed model category, and let
\[
a\colon W\to X, \qquad b\colon X\to Y, \qquad c\colon Y\to Z
\]
be maps in $\EC$.  If $b\circ a$ and $c\circ b$ are weak equivalences,
then so are $a$, $b$, and $c$.
\end{proposition}

\begin{proof}
For any object $V$ in $\EC$, the map $(b\circ a)_{*}\colon
\Ho\EC(V,W)\to \Ho\EC(V,Y)$ is a bijection, and so $b_{*}\colon
\Ho\EC(V,X)\to \Ho\EC(V,Y)$ is a surjection.  The map $(c\circ
b)_{*}\colon \Ho\EC(V,X)\to \Ho\EC(V,Z)$ is a bijection, and so
$b_{*}\colon \Ho\EC(V,X)\to \Ho\EC(V,Y)$ is an injection. Thus,
$b_{*}$ is a bijection for every $V$ in $\EC$, and so by the Yoneda
Lemma, $b$ is an isomorphism in $\Ho\EC$. It follows that $b$ is a
weak equivalence, and by the two out of three axiom, that $a$ and $c$
are weak equivalences.
\end{proof}


\section{Proofs of the main results}\label{PrfSec}

In this section, we apply the theory of enriched derived functors
developed in the Sections~\ref{EnrModCatSec}--\ref{BiFunSec} to prove
the theorems stated in the introduction.  Throughout, we assume that
$\Moncat$ is a monoidal model category.  We use $\monoid$ generally to
denote an arbitrary monoid in $\Moncat$ and $\monoidtwo$ to denote a
monoid in whose underlying object in $\Moncat$ is semicofibrant.
Also, we assume that all of the categories of modules being discussed
are closed model category with fibrations and weak equivalences
created in $\Moncat$.

In order to complete the proofs of the results stated in the
introduction, we must show that our results on enriched parametrized
adjunctions can be applied to the fundamental parametrized adjunctions
arising in the study of bimodules.  The following proposition provides
a general statement.

\begin{proposition}\label{propintmain}
Let $\monoidone$, $\monoidtwo$, and $\monoidthree$ be monoids in
$\Moncat$.  If the underlying object of $\monoidtwo$ is semicofibrant
in $\Moncat$, then:
\begin{thmlist}
\item For cofibrations $f\colon M\to M'$ in
\Bimod{\monoidone}{\monoidtwo} and $g\colon N\to N'$ in
\Bimod{\monoidtwo}{\monoidthree}, the map
\[
(M\balprod[\monoidtwo]N')\cup_{(M\balprod[\monoidtwo]N)}
(M'\balprod[\monoidtwo]N)\to M'\balprod[\monoidtwo] N'
\]
is a cofibration in \Bimod{\monoidone}{\monoidthree} and is a weak
equivalence if either $f$ or $g$ is.
\item For a cofibration $f\colon M\to M'$ in
\Bimod{\monoidone}{\monoidtwo} and a fibration $p\colon P'\to P$ in
\Bimod{\monoidone}{\monoidthree}, the map
\[
\lmfunc{M'}{P'}\to \lmfunc{M'}{P}\times_{\lmfunc{M}{P}}\lmfunc{M}{P'}
\]
is a fibration in \Bimod{\monoidtwo}{\monoidthree} and is a weak
equivalence if either $f$ or $p$ is.

\item For a cofibration $g\colon N\to N'$ in
\Bimod{\monoidtwo}{\monoidthree} and a fibration $p\colon P'\to P$ in
\Bimod{\monoidone}{\monoidthree}, the map
\[
\rmfunc[\monoidthree]{N'}{P'}\to
\rmfunc[\monoidthree]{N'}{P}\times_{\rmfunc[\monoidthree]{N}{P}}
\rmfunc[\monoidthree]{N}{P'}
\]
is a fibration in \Bimod{\monoidone}{\monoidtwo} and is a weak
equivalence if either $g$ or $p$ is.
\end{thmlist}
\end{proposition}

\begin{proof}
By the usual Quillen adjunction argument, part~(a) is equivalent to
both part~(b) and part~(c).  Part~(b) follows from
Proposition~\ref{SemiCoPropB} and the Enrichment Axiom for \Leftmod
(Proposition~\ref{leftmodenr}).
\end{proof}

As previously indicated, Theorem~\ref{intenr} is a special case of
Theorem~\ref{mainenr}, and Proposition~\ref{SemiCoPropA} follows from
the results proved in the previous section.  We now go through the
proofs of the remaining theorems from the introduction:

\begin{proof}[Proof of Theorem~\ref{ForgetAdjProp}] Applying
Theorem~\ref{enradj}, part~(a) is clear from the hypothesis that the
fibrations and weak equivalences are created in \Moncat and part~(b)
is a formal consequence of part~(a) since enriched right adjoints
preserve cotensors.  Parts~(c) and~(d) follow similarly from
Theorem~\ref{enradj} and the definition of semicofibrant.
\end{proof}

\begin{proof}[Proof of Theorem~\ref{intexist}] Parts (a) and (b) of
this theorem follow from Theorem \ref{epqa} and Proposition
\ref{propintmain}.  The claim in part (c) of the theorem that
$M\balprod[\monoidtwo](-)$ and $\lmfunc[\monoidone]{M}{-}$ form a
Quillen adjoint pair is just a reformulation of the hypothesis that
$M$ is semicofibrant in \Leftmod.  Theorem \ref{enradj} provides
enriched derived adjoint functors for this Quillen pair.  The
equivalence of $\WowExt{\monoidone}{\monoidtwo}{\monoidthree}{M}{-}$
with the right derived functor is a consequence of
Theorem~\ref{semigood} and the equivalence of
$\WowTor{\monoidone}{\monoidtwo}{\monoidthree}{M}{-}$ with the left
derived functor follows by the uniqueness of left adjoints.
\end{proof}

\begin{proof}[Proof of Theorem~\ref{intnat}]
Note that for each of the natural transformations whose existence is
asserted by this theorem, there is an obvious corresponding natural
transformation before passage to the homotopy categories.  Applying
Theorem \ref{penrquilcri} and the universal property of the enriched
total left and right derived bifunctors to these known natural
transformations yields the desired natural transformations between the
derived functors.
\end{proof}

\begin{proof}[Proof of Theorem~\ref{intforcomp}]
By hypothesis the forgetful functor preserves fibrant objects and by
Proposition~\ref{SemiCoPropB} it preserves cofibrant objects when the
monoid (whose action is being forgotten) is semicofibrant in
$\Moncat$.  The theorem follows by applying
Proposition~\ref{composederived}:
\begin{thmlist}
\item The right adjoint version, with $\Phi=\emptylmfunc$, $\Psi$ the
forgetful functor $\Bimod{\monoidone}{\monoidtwo}\to
\Leftmod[\monoidone]$, and $\Xi$ the forgetful functor
$\Bimod{\monoidone}{\monoidthree}\to \Leftmod[\monoidone]$ in
part~(a).
\item With $\Phi =\balprod[\monoidtwo]$, $\Psi$ the forgetful functor
functor $\Bimod{\monoidone}{\monoidtwo}\to \Rightmod[\monoidtwo]$, and
$\Xi$ the identity functor in part~(b).
\item With $\Phi =\balprod[\monoidtwo]$, $\Psi$ , the identity
functor, and $\Xi$ the forgetful functor
$\Bimod{\monoidtwo}{\monoidthree}\to \Leftmod[\monoidtwo]$ in
part~(c).
\end{thmlist}\unskip\vskip-\baselineskip
\end{proof}

\begin{proof}[Proof of Theorem~\ref{intassoc}] The statement about
$\WTor$ is a straightforward application of
Proposition~\ref{composederived}, which can be applied inductively to
any association.  The statement about $\WExt$ is adjoint.
\end{proof}


\section{Accommodating non-semicofibrant monoids}\label{NonSemiSec}

In the theorems of the introduction we needed to impose the hypothesis
that certain monoids have semicofibrant underlying objects in
$\Moncat$.  While the results there appear to be the best possible for
an arbitrary monoidal model category, the monoidal model categories
used in practice tend to satisfy even stronger properties which allow
the semicofibrancy hypothesis to be partially dropped. Specifically,
in this section we consider monoidal model categories $\Moncat$ where
all categories of modules are closed model categories with fibrations
and weak equivalences created in $\Moncat$, and satisfy in addition
the following properties:
\begin{hyplist}
\item For any monoid $\monoid$, there exists a monoid $\monoid'$ with
underlying object in $\Moncat$ semicofibrant and a map of monoids
$\monoid'\to\monoid$ that is a weak equivalence.
\item For any monoid $\monoid$ and any cofibrant left $\monoid$-module
$M$, the functor $(-)\balprod[\monoid] M$ preserves weak equivalences
between all right $\monoid$-modules.
\end{hyplist}
\noindent \emph{For the statements in this section, the monoidal model category
$\Moncat$ is always assumed to satisfy properties~{(i)} and~{(ii)} above.}
\smallskip

The first property holds in particular when the conclusions of
\cite[4.1]{ssmonoidal} hold: The category of monoids in $\Moncat$ is
then itself a closed model category and the cofibrant objects have
their underlying object in $\Moncat$ semicofibrant.  Although the we
do not know of a general principle that would imply the second
property, it holds in all presently known monoidal model categories of
spectra \cite{EKMM,HSS,FSPs} and equivariant spectra on complete
universes \cite{LMuct,GOrtho,gss} as well as the
most common monoidal model categories coming from algebra.  The
purpose of this section is to indicate specifically which of the
semicofibrancy hypotheses of the theorems of the introduction can be
eliminated under the assumptions above.

Theorem~\ref{intenr} requires no semicofibrancy hypothesis.
Property~(ii) above, applied with the monoid $\unit$, shows that the
comparison map between tensors in $\Ho\Moncat$ and tensors in
$\Ho\Leftmod$ is a natural isomorphism.

\begin{theorem}
Let $\monoid$ be a monoid in $\Moncat$.  Then $\Ho\Leftmod$ is enriched over
$\Ho\Moncat$ by the right derived functor $\emptyderlmfunc$ of
$\emptylmfunc$, tensored by the left derived functor of $\monprod$ and
cotensored by the right derived functor of $\emptyfunc$.  Moreover,
the derived forgetful functor $\Ho\Leftmod\to\Ho\Moncat$ preserves
tensors and cotensors.
\end{theorem}

Property~(ii) above implies that $\balprod$ satisfies the hypotheses
of Theorem~\ref{penrquilcri}, and we can therefore define
$\WTor_{\monoid}$ to be its enriched total left derived bifunctor.  In general,
$\balprod$ does not satisfy the hypotheses
of Theorem~\ref{epqa}, and we we need to work a bit harder to find a
right adjoint.  Applying property~(i) above
to find a weak equivalence $\monoid'\to\monoid$ with $\monoid'$
semicofibrant in $\Moncat$, property~(ii) implies both that the
extension of scalars and forgetful functor adjunction between
$\Leftmod$ and $\Leftmod[\monoid']$ is a Quillen equivalence and also
that the natural transformation $\balprod[\monoid']\to \balprod$
induces an enriched natural isomorphism of left derived functors
$\WTor_{\monoid'}\to\WTor_{\monoid}$.  This implies that
$\WTor_{\monoid}$ fits into an enriched parametrized adjunction.  The
right adjoint is a refinement of $\WExt_{\unit}$ and so has some
justification to be denoted as $\WowExt{}{\monoid}{}{-}{-}$, but in
general will not be the right derived functor of
\[
\emptyfunc \colon \Rightmod \times \Moncat \to \Leftmod.
\]
Since comparison map
$\WowExt{}{\monoid'}{}{-}{-}\overto{\iso}\emptyderfunc$ is adjoint to
the map $\WowTor{}{}{}{-}{-}\to \WowTor{}{\monoid'}{}{-}{-}$ induced
by $\monprod \to\balprod$, the map
$\WowExt{}{\monoid}{}{-}{-}\overto{\iso}\emptyderfunc$ has an
analogous description. We summarize this in the following theorems.

\begin{theorem}
Let $\monoid$ be a monoid in $\Moncat$.
\begin{thmlist}
\item\label{badtor} The total left derived bifunctor $\WTor_{\monoid}$ of
$\balprod$ exists, is enriched over $\Ho\Moncat$, and is an enriched
parametrized left adjoint in each variable.
\item\label{badext} The right adjoints $\WowExt{}{\monoid}{}{-}{-}$
and $\WowExtR{}{\monoid}{}{-}{-}$ are naturally isomorphic to
$\derfunc{-}{-}$ in $\Ho\Moncat$ by the adjoint to the comparison map
$\WowTor{}{}{}{-}{-}\to \WowTor{}{\monoid}{}{-}{-}$.
\item\label{btorgood}For each fixed right module $M$ and each fixed
left module $N$, $\WowTor{}{\monoid}{}{M}{-}$ and
$\WowTor{}{\monoid}{}{-}{N}$ are the left derived functors of
$M\balprod(-)$ and $(-)\balprod N$.
\end{thmlist}
\end{theorem}

\begin{theorem}
Let $\monoid'\to \monoid$ be a map of monoids and a weak equivalence
in $\Moncat$.  Then the forgetful functor $\Leftmod\to\Leftmod[\monoid']$ is
the right adjoint of a Quillen equivalence.  The derived equivalence
of homotopy categories preserves tensors and cotensors, and the
universal enriched natural transformation $\WTor_{\monoid'}\to
\WTor_{\monoid}$ is a natural isomorphism.
\end{theorem}

The bimodule version of $\WTor$ is complicated by the fact that a pair
of weak equivalences of monoids $\monoidone'\to\monoidone$ and
$\monoidtwo'\to\monoidtwo$ does not necessarily induce a weak
equivalence $\monoidone'\monprod \monoidtwo'\to \monoidone\monprod
\monoidtwo$, and so does not necessarily induce a Quillen equivalence
between categories of bimodules.  However, it follows from
property~(ii) above, that the map is a weak equivalence when one of
$\monoidone'$, $\monoidtwo'$ and one of $\monoidone,\monoidtwo$ are
semicofibrant in $\Moncat$.

\begin{proposition}
Let $\monoidone'\to\monoidone$ and $\monoidtwo'\to\monoidtwo$ be maps
of monoids in $\Moncat$.  Then the forgetful functor
$\Bimod{\monoidone}{\monoidtwo}\to \Bimod{\monoidone'}{\monoidtwo'}$
is the right adjoint of a Quillen adjunction.  If both maps are weak
equivalences and one of $\monoidone'$, $\monoidtwo'$ and one of
$\monoidone,\monoidtwo$ are semicofibrant in $\Moncat$, then the
Quillen adjunction is a Quillen equivalence.
\end{proposition}

When $\monoidthree$ is a monoid whose underlying object is
semicofibrant in $\Moncat$, then cofibrant
$(\monoidtwo,\monoidthree)$-bimodules are cofibrant as left
$\monoidtwo$-modules.  Applying property~(ii) again, we obtain the
following refinement of the theorems from the introduction.

\begin{theorem}\label{mainbad}
Let $\monoidone$, $\monoidtwo$, and $\monoidthree$ be monoids in $\Moncat$ 
and assume that $\monoidthree$ is semicofibrant in $\Moncat$.  Then
the left derived functor
$\WowTor{\monoidone}{\monoidtwo}{\monoidthree}{-}{-}$ exists, is an
enriched parametrized left adjoint in the second variable
(parametrized by the first variable), and the enriched natural map
\[
\WowTor{}{\monoidtwo}{}{-}{-}\to
\WowTor{\monoidone}{\monoidtwo}{\monoidthree}{-}{-}
\]
is an isomorphism.  Moreover, if either $\monoidone$ or $\monoidtwo$
is semicofibrant in $\Moncat$, then
$\WowTor{\monoidone}{\monoidtwo}{\monoidthree}{-}{-}$ is also an
enriched parametrized left adjoint in the first variable.
\end{theorem}

Since by \ref{btorgood}, we have that
$\WowTor{\monoid}{\monoid}{}{\monoid}{-}$ and therefore
$\WowExt{\monoid}{\monoid}{}{\monoid}{-}$ are naturally isomorphic to
the identity functor, the previous theorem can be applied as in the
introduction to the case of $\monoidthree =\unit$ to study the
extension of scalars and coextension of scalars functors.

\begin{corollary}
Let $\monoidone\to\monoidtwo$ be a map of monoids in $\Moncat$.  Then
the derived forgetful functor $\Ho\Leftmod[\monoidtwo]\to
\Ho\Leftmod[\monoidone]$ has both a left and a right adjoint.  The
left adjoint is naturally isomorphic in $\Ho\Moncat$ to
$\WowTor{}{\monoidone}{}{\monoidtwo}{-}$ and the right adjoint is
naturally isomorphic in $\Ho\Moncat$ to
$\derlmfunc[\monoidone]{\monoidtwo}{-}$.
\end{corollary}

The map $\unit\to\monoid$ gives the free and cofree functors on
homotopy categories.  Using the map of monoids $\monoid\iso
\monoid\monprod \unit\to \monoid\monprod\monoid$, the universal
property of the left derived functors and the universal property of
right adjoints induce comparison maps between the free functor and the
functors
\[
\derfree (-)=\WowTor{\monoid\monprod\monoid^{\op}}{}{}{\monoid}{-}
\tand
\dercofree(-)=\WowExt{}{\monoid\monprod\monoid^{\op}}{}{\monoid}{-},
\]
to $\Ho\Bimod{\monoid}{\monoid}$.  Since the comparison maps with
$\WowTor{}{}{}{\monoid}{-}$ and $\derfunc{\monoid}{-}$ are
isomorphisms and the derived forgetful functor reflects isomorphisms,
we obtain the first part of the following theorem.  The isomorphisms
in the second part follow because the derived forgetful functor
$\Ho\Leftmod\to\Ho\Moncat$ takes the comparison maps to the
corresponding ones for the free and cofree functors under the natural
isomorphism from the first part.

\begin{theorem}
The free and cofree functors $\Ho\Moncat \to \Ho\Leftmod$ are enriched
naturally isomorphic to the composite of the functors
$\derfree,\dercofree\colon \Ho\Moncat\to\Ho\Bimod{\monoid}{\monoid}$
and the derived forgetful functor
$\Ho\Bimod{\monoid}{\monoid}\to\Ho\Leftmod$.  Moreover, the canonical
comparison maps
\[
\dertensor{M}{X}\to \WowTor{\monoid}{\monoid}{}{\derfree X}{M}\tand
\dercotensor{X}{M}\to \WowExt{\monoid}{\monoid}{}{\derfree X}{M}
\]
in $\Ho\Leftmod$, and
\[
\WowExt{\monoid}{}{\monoid}{M}{\dercofree X}\to
\WowExtR{}{\monoid}{}{M}{X}
\]
in $\Ho\Rightmod$ are isomorphisms.
\end{theorem}


\bibliographystyle{plain}

\begin{thebibliography}{10}

\bibitem{modcat}
W.~G. Dwyer and J.~Spali{\'n}ski.
\newblock Homotopy theories and model categories.
\newblock In {\em Handbook of algebraic topology}, pages 73--126.
  North-Holland, Amsterdam, 1995.

\bibitem{DHKS}
William~G. Dwyer, Philip~S. Hirschhorn, Daniel~M. Kan, and Jeffrey~H. Smith.
\newblock {\em Homotopy limit functors on model categories and homotopical
  categories}, volume 113 of {\em Mathematical Surveys and Monographs}.
\newblock American Mathematical Society, Providence, RI, 2004.

\bibitem{EKMM}
A.~D. Elmendorf, I.~Kriz, M.~A. Mandell, and J.~P. May.
\newblock {\em Rings, modules, and algebras in stable homotopy theory},
  volume~47 of {\em Mathematical Surveys and Monographs}.
\newblock American Mathematical Society, Providence, RI, 1997.
\newblock With an appendix by M. Cole.

\bibitem{HSS}
M.~Hovey, B.~Shipley, and J.~Smith.
\newblock Symmetric spectra.
\newblock {\em J. Amer. Math. Soc.}, 13(1):149--208, 2000.

\bibitem{HvModCat}
Mark Hovey.
\newblock {\em Model Categories}, volume~63 of {\em Mathematical Surveys and
  Monographs}.
\newblock American Mathematical Society, Providence, RI, 1999.

\bibitem{BasicEnrch}
G.~M. Kelly.
\newblock {\em Basic concepts of enriched category theory}.
\newblock Cambridge University Press, Cambridge, 1982.

\bibitem{LMuct}
L.~G. Lewis, Jr. and M.~A. Mandell.
\newblock Equivariant universal coefficient and {K}\"unneth spectral sequences.
\newblock {\em Proc. London Math. Soc. (3)}, 92(2):505-544, 2006.

\bibitem{GOrtho}
M.~A. Mandell and J.~P. May.
\newblock Equivariant orthogonal spectra and {$S$}-modules.
\newblock {\em Mem. Amer. Math. Soc.}, 159(755):x+108, 2002.

\bibitem{FSPs}
M.~A. Mandell, J.~P. May, S.~Schwede, and B.~Shipley.
\newblock Model categories of diagram spectra.
\newblock {\em Proc. London Math. Soc. (3)}, 82(2):441--512, 2001.

\bibitem{gss}
Michael~A. Mandell.
\newblock Equivariant symmetric spectra.
\newblock In {\em Homotopy theory: relations with algebraic geometry, group
  cohomology, and algebraic $K$-theory}, volume 346 of {\em Contemp. Math.},
  pages 399--452. Amer. Math. Soc., Providence, RI, 2004.

\bibitem{quil}
D.~G. Quillen.
\newblock {\em Homotopical Algebra}, volume~43 of {\em Lecture Notes in
  Mathematics}.
\newblock Springer, Berlin, 1967.

\bibitem{ssmonoidal}
Stefan Schwede and Brooke~E. Shipley.
\newblock Algebras and modules in monoidal model categories.
\newblock {\em Proc. London Math. Soc. (3)}, 80(2):491--511, 2000.

\end{thebibliography}

\end{document}